\documentclass[11pt]{article}

\usepackage[margin=1in]{geometry}

\usepackage{amsmath}
\usepackage{hyperref}
\usepackage{amsfonts}
\usepackage{amssymb}
\usepackage{graphicx}
\usepackage{multicol,float}
\usepackage{color}
\usepackage{multirow}
\usepackage{amsthm}
\usepackage{dsfont}
\usepackage{setspace}
\usepackage{bm}
\usepackage{caption}
\usepackage{subcaption}
\usepackage{comment}

\usepackage{mathtools}
\mathtoolsset{showonlyrefs=true}
\usepackage{epstopdf}
\DeclareMathAlphabet{\mathpzc}{OT1}{pzc}{m}{it}

\usepackage[numbers]{natbib}

\newtheorem{theorem}{Theorem}[section]

\newtheorem{definition}{Definition}[section]

\newtheorem{remark}{Remark}[section]

\numberwithin{equation}{section}

\newcommand \Del {\Delta}
\newcommand \eps {\varepsilon}
\newcommand \gam {\gamma}

\newcommand \xt {\tilde x}
\newcommand{\id}{\mathds{1}}

\newcommand \Vt {\widetilde V}

\newcommand \drm {\mathrm{d}}
\newcommand \e {\mathrm{e}}
\newcommand \Pay {\mathrm{P}}
\newcommand \NT {\mathrm{NT}}

\newcommand \E {\mathbb{E}}
\newcommand \F {\mathbb{F}}
\renewcommand \P {\mathbb{P}}
\newcommand \R {\mathbb{R}}
\newcommand \V {\mathbb{V}}

\newcommand \mA {\mathcal{A}}
\newcommand \mC {\mathcal{C}}
\newcommand \mF {\mathcal{F}}
\newcommand \mM {\mathcal{M}}

\newcommand \pD {\mathpzc{D}}

\newcommand \pd {\mathpzc{d}}

\DeclareMathOperator*{\argsup}{arg\,sup}
\newcommand \noi {\noindent}

\allowdisplaybreaks

\begin{document}

	\title{Equilibrium Strategies for 
		Singular Dividend Control Problems under the Mean-Variance Criterion}

	\author{Jingyi Cao%
		\thanks{Department of Mathematics and Statistics, York University, Canada. Email: jingyic@yorku.ca}
		\and Dongchen Li%
		\thanks{Department of Mathematics and Statistics, York University, Canada. Email: dcli@yorku.ca}
		\and Virginia R. Young%
		\thanks{Department of Mathematics, University of Michigan, USA. Email: vryoung@umich.edu}
		\and Bin Zou%
		\thanks{Corresponding author. Department of Mathematics, University of Connecticut, USA. Email: bin.zou@uconn.edu}
	}
	
	\date{This version: November 10, 2025\\
	Accepted for publication in \emph{SIAM Journal on Control and Optimization}}
	\maketitle

\begin{abstract}
We revisit the optimal dividend problem of de Finetti by adding a variance term to the usual criterion of maximizing the expected discounted dividends paid until ruin, 
in a singular control framework.  Investors do not like variability in their dividend distribution, and the mean-variance (MV) criterion balances the desire for large expected dividend payments with small variability in those payments.  
The resulting MV singular dividend control problem is time-inconsistent, and we follow 
a game-theoretic approach to find a time-consistent equilibrium strategy.  Our main contribution is a new verification theorem for the novel dividend problem, in which the MV criterion is applied to an integral of the control until ruin, a random time that is endogenous to the problem.   
We demonstrate the use of the verification theorem in two cases for which we obtain the equilibrium dividend strategy (semi-)explicitly, and we provide a numerical example to illustrate our results.

\medskip

\noi \textit{MSC2020 codes}: 49J40, 49L12, 49N70, 91A23, 91G50

\medskip

\noi \textit{Keywords}: Optimal divided problem, mean-variance criterion, singular control, time inconsistency, verification lemma

\end{abstract}
	
\section{Introduction}
\label{sec:intro}

The optimal dividend problem is a classic topic in actuarial and financial mathematics that aims to find the best strategy for a company to distribute dividends to its shareholders. In a seminal work, de Finetti \cite{D1957} proposes to maximize the expected discounted dividend payments up to the company's ruin time. This objective balances the trade-off between paying out dividends earlier and retaining earnings to ensure future growth and maintain financial stability, and it is arguably the most popular criterion in the study of optimal dividends. However, as argued in Avanzi \cite{A2009} (p.239), ``{variability in dividend payments is not well received in the markets},'' and de Finetti's criterion does \emph{not} penalize variability. This motivates us to incorporate a variance term to penalize dividend variability and consider a mean-variance (MV) criterion for finding the optimal dividend strategy.

We consider a dividend-paying company and model its surplus by a Brownian motion with drift, the so-called diffusion model in risk theory (see Grandell \cite{G1991}).\footnote{The diffusion model is a popular choice in the optimal dividend problems; see Asmussen and Taksar \cite{AT1997} and Taksar \cite{T2000} for earlier works and Albrecher et al.\ \cite{AAM2022} and Guan and Xu \cite{GX2024} for more recent contributions. In particular, Cohen and Young \cite{CY2021} show that if the company uses the optimal strategy under the diffusion approximation but for the scaled Cram\'er-Lundberg risk model, then doing so is $\varepsilon$-optimal, and they specify the order of $\varepsilon$ relative to the scaling factor.} Let $D = \{D_t\}_{t \ge 0}$ denote the company's dividend strategy, in which $D_t$ is the cumulative amount of dividends paid up to time $t$.
We adopt the singular control framework and do not restrict dividend payments to be absolutely continuous, resulting in a singular control problem. Given a dividend strategy $D$, define $\tau := \tau(D)$ to be the first time when the company's surplus $X$ reaches zero or less, referred to as the \emph{ruin time}; let $Y_t$ denote the total dividends paid between $t$ and $\tau$, discounted at a constant  rate $\rho > 0$, that is, $Y_t = \int_t^\tau \, \e^{-\rho (s-t)}  \drm D_s$.  In the classical setup of de Finetti, the goal is to find an optimal dividend strategy that maximizes $\E_{x,t} (Y_t)$, the conditional expectation of $Y_t$ given the initial surplus $X_{t^-} = x \ge 0$. As motivated above, we propose an MV objective, namely, $J(x, t; D) = \E_{x,t} (Y_t) - \frac{\gamma}{2} \V_{x,t}(Y_t)$,\footnote{While MV preferences are among the most popular criteria in decision making, an alternative choice is the mean-standard deviation (MSD) $\tilde{J}(x, t; D) := \E_{x,t} (Y_t) - \frac{\gamma}{2} \sqrt{ \V_{x,t}(Y_t) }$. Note that MSD preferences satisfy the scale-invariance property (that is, $\tilde{J}(x, t; \alpha D) = \alpha \tilde{J}(x, t; D)$ for all $\alpha \ge 0$), which is desirable in some applications (see, for instance, Bayraktar et al.\ \cite{BZZ2019} for an equilibrium stopping problem under MSD).} 
in which $\gam > 0$ regulates the penalty on the variability in dividend payments and can be interpreted as a risk aversion parameter. Note that the limiting case of $\gam \to 0^+$ reduces to de Finetti's model.

It is well known that dynamic MV optimization problems, such as the above MV dividend problem, are time-inconsistent  (see, for instance, Bj\"ork and Murgoci \cite{BM2010}).\footnote{Let $\{u^*_s|_{\{ x,t \} }\}_{s \ge t}$ denote the ``optimal''  strategy of an optimization problem obtained under the initial condition $X_{t^-} = x \ge 0$ for all feasible $(x, t)$. This dynamic problem is called time-inconsistent if $u_s^*|_{\{ x_1, t_1\} } \neq u_s^*|_{\{x_2, t_2\}}$ holds for some $s > t_2 > t_1$ and feasible $x_1$,  
in which $x_2 = X^*_{t_2^-}|_{\{ x_1, t_1\}}$ is the state process at time $t_2^-$ under the strategy $\{u^*_s|_{\{x_1, t_1\}}\}_{t_1 \le t < t_2}$.}
In this work, we follow the game-theoretic approach to seek an equilibrium dividend strategy (see Definition \ref{def:sing_timeconsist}). To that end, we first develop and prove a verification theorem (Theorem \ref{thm:verif_sing}) that is tailored to our MV dividend problem and differs from those in the literature. Next, we apply this theorem to obtain the equilibrium strategies in two distinctive scenarios, large $\gamma$ and small $\gamma$ (risk aversion). To be precise, when $\gamma$ exceeds a threshold, we show that the equilibrium strategy is to pay out the entire surplus and declare bankruptcy immediately (Theorem \ref{thm:eqm_pV_eq_x}). This result is consistent with intuition because a sufficiently large $\gamma$ imposes a big penalty on the variance of dividend payments $Y_t$, and the strategy of paying out all surplus yields a zero variance. When $\gamma$ is sufficiently small, we show that the equilibrium strategy is a time-independent barrier strategy with a strictly positive barrier $\xt$ (Theorem \ref{thm:onebarrier}); that is, given the initial surplus $X_{t^-} = x$, the company pays out a lump sum dividend of $\max\{ x - \xt, 0 \}$ at time $t$, and thereafter pays dividends so that the resulting surplus is reflected at the barrier $\xt$.

Finding an equilibrium solution of the time-inconsistent MV dividend singular control problem stated above is new to the time-inconsistent control literature. We are aware of only three papers that study time-inconsistent singular control problems: Liang et al.\ \cite{LLY2024}, Liang and Luo \cite{LL2025}, and Dai et al.\ \cite{DJJ2024}. Liang et al.\ \cite{LLY2024} study an optimal reinsurance problem for an insurer, and the singular control is irreversible reinsurance coverage; in their paper, time inconsistency  arises from \emph{non-exponential discounting} in the  objective (see, for instance, Section 5 in Bj\"ork and Murgoci \cite{BM2010}).
Liang and Luo \cite{LL2025} extend the model in Liang et al.\ \cite{LLY2024} to a Stackelberg reinsurance game and assume that both the insurer and reinsurer are endowed with MV preferences.
Dai et al.\ \cite{DJJ2024} solve an MV portfolio optimization problem with proportional transaction costs in a standard Black-Scholes market. Apart from the obvious difference in the optimization problem itself, our paper also differs from those in defining the admissible strategies and, later, equilibrium strategies. We outline the key points below and refer the reader to Remark \ref{rem:eq_def} for a detailed discussion.
Liang et al.\ \cite{LLY2024} and Liang and Luo \cite{LL2025} define admissible strategies by partitioning the feasible region into the continuation and intervention regions; but, both Dai et al.\ \cite{DJJ2024} and this paper define admissible strategies in a more standard way (see Definition \ref{def:adm}). On the other hand, Dai et al.\ \cite{DJJ2024} impose additional $\alpha$-H\"older continuity assumptions (with $\alpha \in (0,1]$) on the spike perturbations and define equilibrium in the order of $\varepsilon^\alpha$; by comparison, we follow the standard first-order ($\varepsilon$) condition as in  Bj\"ork and Murgoci \cite{BM2010}. The key to achieving the standard weak equilibrium is the estimate in \eqref{eq:int_est}, which allows us to bound the error term by $o(\varepsilon)$. 

In the standard MV literature, the objective is in the form of $\E_{x, t}(X_T) - \frac{\gamma}{2} \V_{x, t}(X_T)$, in which $X_T$ is the controlled state process at the terminal time $T$. For instance, both Liang and Luo \cite{LL2025} and Dai et al.\ \cite{DJJ2024} follow this standard setup and assume $T$ is a fixed constant horizon ($X_T$ is replaced by $\ln X_T$ in Dai et al.\ \cite{DJJ2024}); 
the majority of the research on time-inconsistent MV classical control problems also adopts this setup (see, for instance, Bj\"ork et al.\ \cite{BMZ2014}).\footnote{For such a setup, because the MV objective only depends on the terminal state $X_T$ through its distribution $\mu(\cdot)$, 
the original MV problem can be reformulated as a control problem of McKean-Vlasov type (MKV problem) over the (infinite-dimensional) \emph{distribution space}. The new MKV problem is time-consistent in the time-distribution space $(t, \mu)$, and one may apply a McKean–Vlasov version of the dynamic programming approach to characterize the optimal value function $v(t, \mu)$; see Ismail and Pham \cite{IP2019} for a nice implementation of this method to MV portfolio optimization problems. However, it is not straightforward to apply this approach to our MV dividend problem, because the MV criterion is applied to $Y_t$, an integral of controls over an endogenously determined,  random time $\tau$.} 
Landriault et al.\ \cite{LLLY2018} study MV investment problems over a random horizon $T$, but they assume that this random $T$ is \emph{independent} of the state process $X$ and control. 
However, the MV objective in this paper involves $Y_t= \int_t^\tau \, \e^{-\rho (s-t)}  \drm D_s$, which is an \emph{integral} of the dividend payments from the current time to the ruin time $\tau$; note that $\tau$ is \emph{endogenously} dependent on the company's dividend strategy and surplus process, a striking difference from the exogenous random horizon in Landriault et al.\ \cite{LLLY2018}. 
Kronborg and Steffensen \cite{KS2015} apply the MV objective to the  terminal wealth $X_T$ and an integral of intertemporal consumption (a classical control) over a finite horizon $[t, T]$; by comparison, $Y_t$ is an integral of dividends (a singular control) over a random horizon $[t, \tau]$. Because of the ``natural'' boundary at $T$, the equilibrium value function in Kronborg and Steffensen \cite{KS2015}  takes the linear-quadratic form (see Proposition 3.1 therein); this ansatz plays an important role in finding (explicit) solutions. Note that this form of solution is similar to that of standard MV problems involving \emph{only} the terminal wealth (see Bj\"ork and Murgoci \cite{BM2010} and Bj\"ork et al.\ \cite{BMZ2014}). However, we do not have an \emph{a priori} guess for the form of the value function $V$ in this paper; in fact, $V$ will not be of linear-quadratic form globally (see $V(x)$ over $x < \xt$ in Theorem \ref{thm:onebarrier}). Regarding the equilibrium consumption $c^*$ (under constant risk aversion), Kronborg and Steffensen \cite{KS2015} show that it is a bang-bang control and only depends on whether the risk-free rate is greater than the discount rate, but is \emph{independent} of the state process $X$. By comparison, the equilibrium dividend strategy $D^*$ in this paper is of barrier type and explicitly depends on the surplus $X$; note that the same conclusion holds under the classical control framework in Cao et al.\ \cite{CLYZ2025}.
The integral form of $Y_t$ in our MV objective significantly complicates the study and leads to an extended system of Hamilton-Jacobi-Bellman (HJB) equations that is different from the systems in related works (see Liang and Luo \cite{LL2025} and Dai et al.\ \cite{DJJ2024} for singular control and Bj\"ork et al.\ \cite{BMZ2014} and Landriault et al.\ \cite{LLLY2018} for classical control). In particular, we remark that our MV objective is \emph{not} a special case of the general MV framework proposed in Bj\"ork and Murgoci \cite{BM2010} (see their objective in equation (39)). Because of this integral form over an endogenously determined random horizon, the HJB system in this work involves three functions: the (equilibrium) value function $V$, the first moment function $G(x,t) = \E_{x, t}(Y_t^*)$, and the second moment function $H(x,t) = \E_{x, t} \big( (Y_t^*)^2 \big)$; see equations \eqref{eq:HJB_V}-\eqref{eq:HJB_H} in Theorem \ref{thm:verif_sing}. But for the standard MV setup of terminal $X_T$, the extend HJB system only involves $V$ and $G$, but not $H$ (see, for instance, Theorem 3.1 in Liang and Luo \cite{LL2025}).

This paper also contributes to the literature on optimal dividends. 
Although MV preferences are well adopted in portfolio selection problems (see Bj\"ork et al.\ \cite{BMZ2014}, Landriault et al.\ \cite{LLLY2018}, and Dai et al.\ \cite{DJJ2024}), they are rarely used in the study of optimal dividend problems. To the best of our knowledge, this is the first paper that solves a singular dividend control problem under MV preferences. By comparison, Cao et al.\ \cite{CLYZ2025} study an MV dividend problem under the classical control framework,\footnote{The dividend strategy in Cao et al.\ \cite{CLYZ2025} is $D_t = \theta_t \, \drm t$ for some bounded dividend rate $0 \le \theta_t \le m$ for all $t \ge 0$, as in Section 2 of Asmussen and Taksar \cite{AT1997}. If the dividend rate process is further required to be non-decreasing, this is referred to as a ratcheting constraint; see Angoshtari et al.\ \cite{ABY2019}, Albrecher et al.\ \cite{AAM2022}, and Guan and Xu \cite{GX2024}.}
and this paper differs from that one in at least three aspects: (1) the definition of equilibrium strategies (see the last point in Remark \ref{rem:eq_def}), (2) the verification lemma (see Remark \ref{rem:veri}), and (3) the main results (see Remark \ref{rem:comp}). Interestingly, a numerical example in Section \ref{sec:nume} suggests that the barrier equilibrium strategy of the classical control problem in Cao et al.\ \cite{CLYZ2025} converges to its counterpart of the singular control problem in this paper, as the maximum dividend rate goes to infinity (see Figure \ref{fig:rate_db}).
As mentioned earlier, time inconsistency might also arise from non-exponential discounting, and related studies on optimal dividend include Chen et al.\ \cite{CLZ2014, CLZ2018}, Zhu et al.\ \cite{ZSY2020}, Zhou and Jin \cite{ZJ2020}, and Christensen and Lindensj\"o \cite{CL2022}, among many others. Please see Albrecher and Thonhauser \cite{AT2009} and Avanzi \cite{A2009} for an overview of the research questions on optimal dividend problems.

The rest of this paper is organized as follows. Section \ref{sec:model} presents the model and main problem. We develop and prove the verification theorem in Section \ref{sec:ver} and apply it to obtain equilibrium dividend strategies in Section \ref{sec:eqm}. We conduct a numerical analysis in Section \ref{sec:nume}.
Finally, Section \ref{sec:con} concludes the study.

\section{Model}\label{sec:model}

We fix a filtered probability space $(\Omega, \mF, \F = (\mF_t)_{t \ge 0}, \P)$, in which the filtration $\F$ is generated by a standard one-dimensional Brownian motion $B = (B_t)_{t \ge 0}$. We consider a company that pays dividends to its shareholders and let $D_t$ denote the \emph{cumulative} amount of dividends paid up to time $t$; we call $D = \{ D_t \}_{t \ge 0}$ a {\it dividend strategy}.  We model the company's uncontrolled surplus by a drifted Brownian motion (see, for instance, Asmussen and Taksar \cite{AT1997}). As such, given a dividend strategy $D$, the company's controlled surplus $X = (X_t)_{t \ge 0}$ follows the dynamics 
\begin{align}
	\label{eq:X}
	\drm X_t = a \, \drm t + b \, \drm B_t - \drm D_t, 
\end{align}
in which $a$ and $b$ are positive constants, with $X_0 >0$. Define the company's ruin time by $\tau := \inf \{ t \ge 0: X_t \le 0 \}$. Let $Y_t$ denote the total dividends paid between time $t$ and ruin time $\tau$ under strategy $D$, discounted at a constant rate $\rho > 0$, that is,\footnote{Throughout this paper, all integrals include the possible jumps at the left end point; for example, $Y_t$ in \eqref{eq:Y} equals $\Delta D_t + \int_{t^+}^\tau \, \e^{-\rho(s - t)} \, \drm D_s$. }
\begin{equation}\label{eq:Y}
	Y_t = \int_t^\tau \e^{-\rho(s - t)} \, \drm D_s, \qquad 0 \le t < \tau.
\end{equation}
We set $Y_t = 0$ for all $t \ge \tau$. 
It is obvious that $X$, $\tau$, and $Y_t$ all depend on the company's dividend strategy $D$, and a more precision notation is to write $X^D$, $\tau^D$, and $Y_t^D$,  but we often suppress this dependence for notational simplicity.

Following the literature on time-inconsistent control problems (see, for instance, Bj\"ork and Murgoci \cite{BM2010} and Bj\"ork et al.\ \cite{BMZ2014} on regular controls and Dai et al.\ \cite{DJJ2024} on singular controls), we focus on (Markov) feedback controls in the form of $D_t = \pD(X_{t^-}, t)$ for some deterministic function $\pD$.  We define admissible dividend strategies below.

\begin{definition}
	\label{def:adm}
	A dividend strategy $D = (D_t)_{t \ge 0}$ is called {\rm admissible} if $(1)$ there exists a Borel-measurable, deterministic function $\pD: \R_+^2 \to \R_+$ such that $D_t = \pD(X_{t^-}, t)$, in which $X$ satisfies \eqref{eq:X} under the strategy $D$; $(2)$ $D$ is non-decreasing over time; $(3)$ $\Delta D_t := D_t - D_{t^-} \le X_{t^-}$ $($that is, the company cannot pay more in dividends that it owns$)$; $(4)$ $D_t = D_\tau$ for all $t \ge \tau$ $($that is, there are no dividend payments after ruin$)$; and $(5)$ $Y_t$ in \eqref{eq:Y} is square integrable for all $t \ge 0$. 
\end{definition}

Note that the above definition trivially generalizes from a starting time of $0$ to an arbitrary starting time $t \ge 0$ (assuming $t < \tau$); let $\mA_t$ denote the set of all admissible strategies $D=(D_s)_{s \ge t}$ for every $t \ge 0$. 
With a slight abuse of notation, we use $D$ to denote both the deterministic function $\pD$ and the dividend strategy induced by it via $D_t = \pD(X_{t^-}, t)$.

As argued in Section \ref{sec:intro}, we assume that the manager of the company penalizes variability in dividend payments by their variance and applies the MV criterion when choosing the company's dividend strategy. In particular, the manager's (dynamic) objective function is given by 
\begin{align}
	\label{eq:J}
	J(x, t; D) = \E_{x, t}(Y_t) - \dfrac{\gam}{2} \, \V_{x, t}(Y_t), \quad D \in \mA_t,
\end{align} 
in which $\gam >0$ is the risk aversion parameter. In \eqref{eq:J}, $\E_{x, t}$ and $\V_{x, t}$ denote expectation and variance, respectively, conditional on $X_{t^-} = x \ge 0$, that is, before any possible lump-sum dividend payments at time $t$. If we set $\gam = 0$, then the objective $J$ in \eqref{eq:J} reduces to the one proposed by de Finetti \cite{D1957} and used in many follow-up works (see Albrecher and Thonhauser \cite{AT2009}).

Because of the variance term in \eqref{eq:J}, maximizing $J(x, t; D)$ for all $(x, t) \in \R_+$ leads to a time-inconsistent control problem. We follow an intrapersonal game approach, as in Bj\"ork and Murgoci \cite{BM2010}, and seek a time-consistent equilibrium dividend strategy $D^*$. 
The definition of an equilibrium strategy under a singular control framework is different from its counterpart under a classical (or regular) control framework in Bj\"ork and Murgoci \cite{BM2010}. Below, we formally define $D^*$, and it is similar to the definition of equilibrium in Dai et al.\ \cite{DJJ2024}; for a different definition, see Liang et al.\ \cite{LLY2024} and Liang and Luo \cite{LL2025}.

\begin{definition}\label{def:sing_timeconsist}
	
Fix an arbitrary initial time $t \ge 0$ and an initial surplus $X_{t-} = x > 0$ and assume that ruin has not occurred by time $t$. 	Let $D^* = (D^*_s)_{s \ge t} \in \mA_t$ be an admissible dividend strategy and denote its associated surplus, ruin time, and discounted dividend payments by $X^* := X^{D^*}$, $\tau^* := \tau^{D^*}$, and $Y_t^* := Y_t^{D^*}$, respectively. 
For a positive number $\eps$, a non-negative number $d \in [0, x]$, and a non-decreasing,  continuous function $\delta$ $($of time only$)$ satisfying $\delta(t+\eps) - \delta(t) = O(\eps)$ as $\eps \to 0$, define a perturbed strategy $D^\eps = (D^\eps_s)_{s \ge t}$ by  
\begin{equation}\label{eq:pD_eps}
D^\eps_s =
\begin{cases}
D^*_{t^-} + d + \int_{t^+}^{s \wedge \tau} \, \drm \delta(u), &\quad t \le s < (t + \eps) \wedge \tau, \vspace{0.5em} \\
D^\eps_{(t + \eps)^-} + \int_{t + \eps}^s \drm D^*_u, &\quad s \ge (t + \eps) \wedge \tau.
\end{cases}
\end{equation}
in which $\tau := \tau^{D^\eps}$ is the ruin time under the perturbed strategy $D^\eps$, and $\int_{t + \eps}^s \drm D^*_u = \Delta D^*_{t+\eps} + \int_{(t + \eps)^+}^s \drm D^*_u$.
The strategy $D^*$ is said to be a {\rm time-consistent equilibrium dividend strategy} if, for all $(x, t) \in \R_+^2$, $d \in [0,x]$, and $\delta$ functions that satisfy the above conditions, 
\begin{equation}\label{eq:pJ_eqm}
\liminf_{\eps \to 0^+} \, \dfrac{J(x, t; D^*) - J(x, t; D^\eps)}{\eps} \ge 0,
\end{equation}
and the {\rm equilibrium value function} $V$ equals
\begin{equation}\label{eq:V}
V(x, t) = J(x, t; D^*).
\end{equation}
\end{definition}

We end this section with a technical remark on the definition of the equilibrium strategies $D^*$ above and a discussion on the existence, (non)uniqueness, and ``optimality'' of equilibria. 

\begin{remark}
	\label{rem:eq_def}
	The definition of $D^\eps$ in \eqref{eq:pD_eps} is largely inspired by Dai et al.\ {\rm \cite{DJJ2024}} who also study a time-inconsistent singular control problem, and it shares the same idea as in Bj\"ork and Murgoci \cite{BM2010} under the regular control framework. 
	In Definition {\rm \ref{def:sing_timeconsist}}, we only require $\delta(t+\eps) - \delta(t) = O(\eps)$, and the denominator in \eqref{eq:pJ_eqm} is $\eps$, the first order of the error $\eps$. However, Dai et al.\ {\rm \cite{DJJ2024}} assume that $\delta$ is $\alpha$-H\"older continuous for some $\alpha \in (0,1]$, and the corresponding denominator is $\eps^\alpha$. To our understanding, the ``small'' terms in Dai et al.\ {\rm \cite{DJJ2024}} are \emph{not} of order $o(\eps)$, and that is why they impose the additional assumption of $\alpha$-H\"older continuity and change the denominator from $\eps$ to $\eps^\alpha$ $($see Definition 2 therein$)$.  We can relax their assumption because after carefully collecting all the integral terms of $\delta$ with order $O(\eps)$,  the summation is of a definite sign $($``negative'' in \eqref{eq:int_est}$)$, which allows us to prove the inequality in \eqref{eq:pJ_eqm}. 
	  Although we assume that $\delta$ is a deterministic, univariate function of time only, we can easily generalize to allowing perturbations in the form of $\delta_s := \delta(X_s, s)$ for some bivariate function $\delta$, as long as $\delta(X_{t+\eps}, t+\eps) - \delta(x-d, t) = O(\eps)$ holds uniformly.  Under that extension, the class of $\delta$ would be large enough to incorporate the $($bounded$)$ \emph{dividend-rate} case. Indeed, let $\theta(X_s, s) \in [0, m]$ be the dividend rate paid at time $s \in (t, t+ \eps)$; then, $\int_{t^+}^{t + \eps} \drm \delta(X_u, u) = \int_{t^+}^{t + \eps} \e^{-\rho(s - t)} \theta(X_s, s) \, \drm s = O(\eps)$, and we can easily choose $\theta_s$ so that $\delta$ is \emph{not} $\alpha$-H\"older continuous, versus the requirement in Dai et al.\ {\rm \cite{DJJ2024}}.
	
Careful readers will notice that the only singular perturbation over $[t, t+ \eps)$ occurs at time $t$ in the definition of $D^\eps$ in \eqref{eq:pD_eps}. It is straightforward to extend from one jump at time $t$ to a countable number of jumps over $[t, t+ \eps)$, but this requires the additional assumption of $\sum_{s \in (t, t+\eps)} \, \Delta D_s^\eps = o(\eps)$. Note that Liang et al.\ {\rm \cite{LLY2024}} impose exactly the same assumption in their definition $($see Definition $2.2(c)$, p. $3217)$.\footnote{The assumption of $\sum_{s \in (t, t+\eps)} \, \Delta D_s^\eps = o(\eps)$ seems to be required in Liang and Luo \cite{LL2025} as well, even though they write $O(\eps)$ instead of $o(\eps)$ (see Definition 2.5(c), p.172).}   
	
By Definition {\rm \ref{def:adm}}, if ruin has occurred before time $t$ or $X_{t-} \le 0$, we have $J(x, t; D) = 0$ for all $D \in \mA_t$. Therefore, to avoid such trivial scenarios, we assume, without loss of generality, that $x > 0$ and ruin has not occurred by time $t$ in Definition {\rm \ref{def:sing_timeconsist}}.   

Recall that Cao et al.\ {\rm \cite{CLYZ2025}} adopt the classical control framework and require admissible dividend strategies to be absolutely continuous. As such, the perturbed strategy $D^\eps$ therein does \emph{not} allow singular jumps, which is equivalent to setting $d \equiv 0$ in \eqref{eq:pD_eps}. In addition, they assume a \emph{linear} form for $\delta$ functions to define perturbed strategies $D_s^\eps$, which, under our notation, yields $D_s^\eps = D_t^* + \int_t^\eps \,  c \, \drm u$ for an arbitrary positive constant $c$ $($less than the maximum dividend rate$)$. Apparently, the perturbed strategies considered in Cao et al.\ {\rm \cite{CLYZ2025}} are  special cases of \eqref{eq:pD_eps}, which allows not only  singular jumps $d > 0$ but also general forms for $\delta$ functions.  

Our definition of the equilibrium strategy in \eqref{eq:pJ_eqm} is the so-called \emph{weak equilibrium}, and it is inspired by the popular approach introduced in Bj\"ork and Murgoci {\rm \cite{BM2010}}. However, one potential drawback of such an approach is that the first-order condition $($FOC$)$ in \eqref{eq:pJ_eqm} is only a necessary condition to characterize equilibrium, and if the FOC holds with equality, there might exist counterexamples in which $J(x, t; D^\eps) - J(x, t; D^*) > 0$ for some small $\eps$, contracting the concept of equilibrium. To address this issue, different notions of equilibrium have been proposed in the literature; see Huang and Zhou {\rm \cite{HZ2021}} and He and Jiang {\rm \cite{HJ2021}} for time-inconsistent control problems, and Bayraktar et al.\ {\rm \cite{BZZ2021}} and Bayraktar et al.\ {\rm \cite{BWZ2023}} for time-inconsistent stopping problems. In this paper, we choose the notion of weak equilibrium because it requires minimal assumptions on the model, and for MV problems, weak equilibria can be characterized by the extended HJB equations $($see Theorem {\rm \ref{thm:verif_sing}} below$)$. By comparison, stronger notions of equilibria require more restrictive assumptions on the model, and they may fail to exist $($see, for instance, Section 4.4 in He and Jiang {\rm \cite{HZ2021}}$)$. Thus, a weak equilibrium is often the first choice when studying a time-inconsistent control or stopping problem, and one proceeds to stronger notions only when there is a good understanding of weak equilibria. As mentioned in the Introduction, the research on time-inconsistent singular control problems is in its early stage, and it is, thus, not surprising that several recent papers $($see Liang et al.\ {\rm \cite{LLY2024}}, Liang and Luo {\rm \cite{LL2025}}, Dai et al.\ {\rm \cite{DJJ2024}}, and Cao et al.\ {\rm \cite{CLYZ2025}}$)$ all choose the notion of weak equilibrium. 
\qed
\end{remark}

\begin{remark}
	As nicely noted in Bj\"ork and Murgoci {\rm \cite{BM2010}}, for all time-inconsistent control problems over a finite, discrete-time horizon, equilibrium strategies, defined similar to the one in Definition {\rm \ref{def:adm}}, always exist and can be obtained by backward recursion, which in turn implies the uniqueness of the equilibrium value function $V$ $($although there may exist multiple equilibrium strategies achieving the same $V)$. However, the existence result is highly nontrivial for the infinite horizon case, due to the lack of natural boundaries, which is shared by our random horizon setup. In addition, uniqueness on $V$ may fail on infinite horizon time-inconsistent problems, and there are concrete examples in the literature that admit multiple equilibria. For instance, Example $3.1$ in Landriault et al.\ {\rm \cite{LLLY2018}} shows that there exist multiple linear equilibrium strategies, each yielding a different $V$, for their MV investment problems over an exponentially distributed random horizon. For the same reason, there is no guarantee on the uniqueness of the equilibrium value function $V$ defined by \eqref{eq:V}. \qed
\end{remark}

\section{Verification theorem}\label{sec:ver}

In this section, we prove a verification theorem for the equilibrium value function $V$ in \eqref{eq:V} and the corresponding equilibrium strategy $D^*$.  We define a differential operator $\mM$ by
\begin{align}\label{eq:mD}
\mM \phi(x, t) = \partial_t \phi(x, t) + a \, \partial_x \phi(x, t) + \dfrac{1}{2} \, b^2 \partial_{xx} \phi(x, t),
\end{align}
in which $\phi \in \mC^{2,1}(\R_+^2)$ and $\partial_{\cdot} \phi$ denotes the corresponding partial derivative of $\phi$.  Because the following verification theorem is relatively new in the literature, we provide its proof in full detail.

\begin{theorem}\label{thm:verif_sing}
Let $\Vt$, $G$, and $H$ be three functions, all mapping from $(x, t) \in \R_+^2$ to $\R$. Define the {\rm pay region} $\Pay$ and {\rm no-transaction region} $\NT$, respectively, by
\begin{equation}\label{eq:P}
	\Pay = \big\{ (x, t) \in \R_+^2: \partial_x \Vt(x, t) = 1 \big\} \quad \text{ and } \quad \NT = \R_+^2 \backslash \Pay.
\end{equation}
Suppose that $\Vt$, $G$, and $H$ satisfy the following conditions: 

\begin{itemize}
\item[$1.$] $\Vt$, $G$, and $H \in \mC^{2,1}(\R_+^2)$, except that $G(\cdot, t)$ and $H(\cdot, t)$ might only be $\mC^1$ along a specific path $x = \xt(t)$ for all $t \ge 0$, with both left and right second derivatives. 

\item[$2.$] $G$ and $H$ satisfy regularity conditions such that the stochastic integrals in \eqref{eq:pG_rando} and \eqref{eq:hat_H} have zero $($conditional$)$ expectation and $\lim_{s \to \infty} \, \E_{x, t} \big(\e^{-\rho (s -t)} \, \phi(X_s, s) \big) = 0$ for $\phi = G, H$.

\item[$3.$] For all $(x, t) \in \R_+^2$, $\Vt$, $G$, and $H$ jointly solve the extended HJB system: 
\begin{align}
		\max \left\{\mM \Vt - \dfrac{\gam}{2}  \, \mM G^2 + \gam G \cdot \mM G - \rho G + \gam \rho \big(H - G^2\big), \; 1 - \partial_x \Vt \right\} &= 0,  \label{eq:HJB_V}
		\\
		\big(\mM G(x, t) - \rho G(x, t) \big) \id_{\{(x, t) \in \NT\}} + \big( 1 - \partial_x G(x, t) \big) \id_{\{(x, t) \in \Pay \}} &= 0, \label{eq:HJB_G}
		\\
		\big(\mM H(x, t) - 2 \rho H(x, t) \big) \id_{\{(x, t) \in \NT\}} + \big( 2G(x, t) - \partial_x H(x, t) \big) \id_{\{(x, t) \in \Pay \}} &= 0, \label{eq:HJB_H}
\end{align}
in which the argument $(x,t)$ is suppressed in \eqref{eq:HJB_V}, with the boundary conditions 
\begin{align}
	\label{eq:boundary}
	\Vt(0, t) = G(0, t) = H(0, t) = 0, \quad \text{for all } t \ge 0.
\end{align}
In addition, there exists an admissible dividend strategy $D^* = (D_s^*)_{s \ge t}$ that solves the Skorokhod reflection problem
\begin{align}
	\label{eq:Sk}
	\begin{cases}
		\drm X_s^* = a \, \drm s + b \, \drm B_s - \drm D_s^*,  & \text{with } X_{t^-}^* = x, \\
		(X_s^*, s) \in \overline{\NT}, & \\
		D_s^* = D_{t^-}^* + \int_t^s \, \id_{\{(x, t) \in \Pay \}} \, \drm D_u^*, &
	\end{cases}
\end{align}
for all $s \ge t$, in which $\overline{\NT}$ denotes the closure of $\NT$ in \eqref{eq:P}. 
\end{itemize}
Then, $\Vt$ is an equilibrium value function defined in \eqref{eq:V}, and $D^*$ is a time-consistent equilibrium dividend strategy. Moreover, $G$ and $H$ have the representations
\begin{align}
	\label{eq:GH}
	G(x, t) = \E_{x, t} \left(Y_t^* \right) \quad \text{and} \quad  H(x, t) = \E_{x, t} \big( \left(Y_t^* \right)^2 \big), 
\end{align}
in which $Y^*$ is the discounted dividends under $D^*$; thus, $V(x,t) = G(x,t) - \frac{\gam}{2} \left( H(x,t) - G^2(x,t) \right)$.
\end{theorem}

Before we prove Theorem \ref{thm:verif_sing}, we provide some intuition for the results. Assume that a lump-sum dividend payment is optimal at $(x, t)$; then, the amount to be paid equals $\argsup_{d \ge 0} \, V(x - d, t) + d$, which motivates the definition of the ``pay region'' $\Pay$ in \eqref{eq:P}. For the ``no-transaction region'' $\NT$, the value function satisfies a standard differential equation, namely, $\mM V - \frac{\gam}{2}   \mM G^2 + \gam G \cdot \mM G - \rho G + \gam \rho \big(H - G^2\big) = 0$. Together, they explain the variational inequality in \eqref{eq:HJB_V} satisfied by $V$. Similarly, both $G$ and $H$ are characterized separately for $(x, t) \in \Pay$ and $(x, t) \in \NT$, leading to \eqref{eq:HJB_G} and \eqref{eq:HJB_H}, respectively.   
Based on the partition of $\Pay$ and $\NT$, we know that if $(x, t) \in \Pay^o$ (interior of $\Pay$), then the manager should immediately pay dividends to reach the boundary of the ``no-transaction region'' $\NT$ or pay out all of $x$ in dividends if $\partial \NT$ (the boundary of $\NT$) is unreachable.\footnote{If $\partial \NT$ is unreachable from $\Pay$, which could occur if $(0, x] \subset \Pay$ for some $x > 0$, then $\Del D^*_t = X_{t^-} = x$, and ruin occurs immediately.}
Thereafter, the interventions are of ``local-time type,'' as described by the third equation in \eqref{eq:Sk}, to keep the company's surplus within the no-transaction region (that is, $(X_s^*, s) \in \overline{\NT}$). Similar to Liang and Luo \cite{LL2025}, the state-time space is divided into two regions (see \eqref{eq:P} here and (3.1)--(3.2) therein), but Dai et al.\ \cite{DJJ2024} further separate the pay region $\Pay$ into ``buy'' and ``sell'' regions in their transaction costs model because buying and selling the risky asset incur costs at different rates.

\begin{proof}
Suppose that $\Vt$, $G$, and $H$ satisfy the conditions of this theorem, and suppose there exists a solution $D^* \in \mA_t$ to the Skorokhod reflection problem in \eqref{eq:Sk}. We prove the theorem in four steps. 

\medskip 
\noi 
\textbf{Step 1.} We show that if $G$ solves \eqref{eq:HJB_G} with $G(0,t) = 0$, then $G(x, t) = \E_{x, t}(Y_t^*) $ in \eqref{eq:GH}. 

Fix  $(x, t) \in \R_+^2$ and a positive number $k > t$.  
By applying It\^o's formula to $\e^{-\rho(\cdot - t)} G(X^*_{\cdot}, \cdot)$, we obtain
\begin{align}\label{eq:pG_rando}
&\e^{-\rho((\tau^* \wedge k) - t)} G(X^*_{\tau^* \wedge k}, \tau^* \wedge k) = G(x, t) + \int_t^{\tau^* \wedge k} \e^{-\rho (s - t)} \big( \mM G(X^*_s, s) - \rho G(X^*_s, s) \big) \, \drm s \notag \\
&\quad   -  \int_t^{\tau^* \wedge k} \e^{-\rho (s - t)} \partial_x G(X^*_s, s) \, \drm D^{*,c}_s + \int_t^{\tau^* \wedge k} \e^{-\rho (s - t)} b \, \partial_x G(X^*_s, s) \, \drm B_s  \notag \\
&\quad   + 
\sum_{s \in [t, \tau^* \wedge k]} \e^{-\rho (s - t)} \big(G(X^*_{s^-} - \Delta D_s^*, s ) - G(X^*_{s^-}, s )\big),
\end{align}
in which $D^{*,c}$ is the continuous part of $D^*$.  The first integral in \eqref{eq:pG_rando} equals $0$ because $(X^*_s, s) \in \overline{\NT}$ for all $s > t$, on which $\mM G - \rho G = 0$ by \eqref{eq:HJB_G}. 
The above discussion implies that a lump-sum dividend ($\Del D^*_s >0$) is only possible at the initial time $t$ when $(x, t) \in \Pay^o$;\footnote{The subsequent analysis follows even without explicitly using this result (that is, temporarily allowing $\Delta D_s^* >0$ for $s >t$). In that case, note $D_s^* = D^{*,c}_s + \sum_{u \in [t, s]} \, \Delta D^*_u$.} 
in that case, we have $G(x - \Del D_t^*, t) - G(x, t) = \int_0^{\Del D_t^*} \, \partial_x G(X^*_{t^-} - u, t) \cdot  \id_{\{(X^*_{t^-}, t ) \in \Pay^o\}} \drm u = \Del D_t^* \cdot  \id_{\{(X^*_{t^-}, t ) \in \Pay^o\}}$ because $\partial_x G = 1$ on $\Pay^o \subset \Pay$. By using this result, we get 
\begin{align*}
	&-  \int_t^{\tau^* \wedge k} \e^{-\rho (s - t)} \partial_x G(X^*_s, s) \, \drm D^{*,c}_s + \sum_{s \in [t, \tau^* \wedge k]} \e^{-\rho (s - t)} \big(G(X^*_{s^-} - \Delta D_s^*, s ) - G(X^*_{s^-}, s )\big) \\
	&= - \int_t^{\tau^* \wedge k} \e^{-\rho (s - t)} \partial_x G(X^*_s, s)  \, \drm D^{*,c}_s - \Del D_t^* = - \int_t^{\tau^* \wedge k} \e^{-\rho (s - t)} \, \drm D^{*}_s,
\end{align*}
in which the last equality uses $\partial_x G = 1$ on $\{ \drm D^{*,c}_s > 0\} \subset \Pay$ and $D^*_s = D^{*,c}_s + \Delta D^*_t$. 

Next, we take conditional expectation on both sides of \eqref{eq:pG_rando} and use the above results, Condition 2 in the theorem, and $G(0, t) =0$ to obtain 
\begin{align*}
	G(x, t) &= \E_{x, t} \bigg(  \int_t^{\tau^* \wedge k} \e^{-\rho (s - t)} \drm D^*_s \bigg) +  \e^{-\rho(k - t)} G(X^*_{k}, k) \cdot \id_{\{ \tau^* > k \}},
\end{align*}
which yields the desired assertion by sending $k \to \infty$ and using Condition 2.

\medskip
\noi 
\textbf{Step 2.} We show that if $H$ solves \eqref{eq:HJB_H} with $H(0,t) = 0$, then $H(x, t) = \E_{x, t} \big( (Y_t^*)^2 \big)$ in \eqref{eq:GH} holds.

Fix $(x, t) \in \R_+^2$ and assume that ruin has not occurred by time $t$. Define a sequence of stopping times $\{\eta_n\}_{n=1,2,\dots}$ by $\eta_n := \inf \{ s \ge t : X_s^* \ge n \}$. For a fixed $k > t$, denote $\tau_{n, k} = \tau^* \wedge k \wedge \eta_n$; define functions $\hat G$ and  $\hat H$ by
\[
\hat G(x, t) = \e^{-\rho t} G(x, t), \quad \hbox{and} \quad 
\hat H(x, t) = \e^{-2\rho t} H(x, t).
\]
By applying It\^o's formula to $\hat G(X^*_{\cdot}, \cdot)$ as in Step 1, we deduce
\begin{equation}\label{eq:hatpG_rando}
\hat G(X^*_s, s) =  \int_s^{\tau_{n,k}} \e^{-\rho u} \, \drm D^*_u -  \int_s^{\tau_{n,k}}  b \, \hat \partial_x G(X^*_u, u) \, \drm B_u + \hat G(X^*_{\tau_{n,k}}, \tau_{n,k}).
\end{equation}
It follows from \eqref{eq:HJB_H} and \eqref{eq:boundary} that $\hat H$ solves
\begin{equation}\label{eq:hatpH}
\mM \hat H(x, t) \id_{\{(x, t) \in \NT\}} + \big( 2\e^{-\rho t} \hat G(x, t) - \partial_x \hat H(x, t) \big) \id_{\{(x, t) \in \Pay\}} = 0,  \quad 
\hat H(0, t) = 0.
\end{equation}
By applying It\^o's formula to $\hat H (X^*_{\cdot}, \cdot )$ and using the results from Step 1, we obtain 
\begin{align}\label{eq:hat_H}
\hat H(X^*_{\tau_{n,k}}, \tau_{n,k}) &= \hat H(x, t) + \int_t^{\tau_{n,k}} \mM \hat H (X^*_s, s) \, \drm s +  \int_t^{\tau_{n,k}}  b \, \partial_x \hat{H} (X^*_s, s) \, \drm B_s \notag \\
&\quad -  \int_t^{\tau_{n,k}} \partial_x \hat{H} (X^*_s, s) \, \drm D^{*,c}_s + \sum_{s \in [t , \tau_{n,k}]} \left( \hat{H} (X^*_{s^-} - \Delta D^*_s, s) - \hat{H} (X^*_{s^-} , s) \right) \notag \\
&= \hat H(x, t) + \int_t^{\tau_{n,k}} \mM \hat H (X^*_s, s) \, \drm s +  \int_t^{\tau_{n,k}}  b \, \partial_x \hat{H} (X^*_s, s) \, \drm B_s \notag \\
&\quad  -  \int_t^{\tau_{n,k}} \partial_x \hat{H} (X^*_s, s) \, \drm D^{*}_s + \left(\partial_x \hat{H} (x, t) \Delta D^*_t + \hat{H} (x - \Delta D^*_t, t) - \hat{H} (x , t) \right)  \id_{\{(x, t) \in \Pay^o\}} \notag \\
&= \hat H(x, t) +  \int_t^{\tau_{n,k}}  b \, \partial_x \hat{H} (X^*_s, s) \, \drm B_s -  \int_t^{\tau_{n,k}} 2\e^{-\rho s} \hat{G} (X^*_s, s) \, \drm D^{*}_s  \notag \\
&\quad  + \left(2\e^{-\rho t} \hat{G}(x, t) \Delta D^*_t + \hat{H}(x - \Delta D^*_t, t) - \hat{H}(x, t) \right)  \id_{\{(x, t) \in \Pay^o\}},
\end{align}
in which we use $\mM \hat{H} = 0$ on $\NT$ and $\partial_x \hat{H} = 2\e^{-\rho t} \hat{G}$ on $\Pay$; recall that $X^*_{s^-} - X^*_s = \Delta D^*_s > 0$ if and only if $s = t$ and $(x, t) \in \Pay^o$. 
To analyze $\hat{H}(x - \Delta D^*_t, t) - \hat{H}(x, t)$ when $(x, t) \in \Pay^o$, note that for all $z \in [x - \Delta D^*_t, x]$,  $(z, t) \in \Pay$, and given $X_{t^-} = z$, there is an immediate lump-sum payment of size $z - (x - \Delta D^*_t)$ at time $t$, implying 
$ G(z, t) = \E_{z,t} (Y_t^*) =  \E_{x - \Delta D^*_t, t}(Y^*_t) + \big(z - (x - \Delta D^*_t) \big)$. 
Using these results, along with \eqref{eq:hatpH},  we have 
\begin{align}\label{eq:DelpH}
\hat{H}(x - \Delta D^*_t, t) - \hat{H}(x, t) &= - \int_{x - \Delta D^*_t}^x \partial_x \hat{H} (z, t) \, \drm z = - 2 \e^{-\rho t} \int_{x - \Delta D^*_t}^x \hat{G}(z, t) \, \drm z \notag \\
&= - 2 \e^{-2\rho t} \int_{x - \Delta D^*_t}^x \left( \E_{x - \Delta D^*_t, t}(Y^*_t) + \big(z - (x - \Delta D^*_t) \big) \right) \drm z \notag \\
&= - 2 \e^{-2\rho t} \left( \Delta D^*_t \cdot \E_{x - \Delta D^*_t, t}(Y^*_t) + \dfrac{1}{2} \, (\Delta D^*_t)^2 \right) \notag \\
&= - \e^{-2\rho t} \left( 2\Delta D^*_t \cdot \E_{x, t} \bigg(\int_t^{\tau^*} \e^{-\rho(s - t)} \, \drm D^*_s  - \Delta D^*_t \bigg) + (\Delta D^*_t)^2 \right) \notag \\
&= -  \e^{-2\rho t} \left( 2\Delta D^*_t \cdot \E_{x, t} \bigg(\int_t^{\tau^*} \e^{-\rho(s - t)} \, \drm D^*_s \bigg) - (\Delta D^*_t)^2 \right).  \hspace{5em}
\end{align}
Then, combining \eqref{eq:hat_H} with \eqref{eq:hatpG_rando} and \eqref{eq:DelpH} and taking conditional expectations imply
\begin{align*}
\hat{H}(x, t) &= \E_{x, t} \Big( \hat{H}(X^*_{\tau_{n,k}}, \tau_{n,k}) \Big)  + 2 \, \E_{x, t} \bigg( \int_t^{\tau_{n,k}} \e^{-\rho s} \left(  \int_s^{\tau_{n,k}} \e^{-\rho u} \drm D^*_u \right) \drm D^*_s \bigg)  \\
&\quad  + 2 \, \E_{x, t} \bigg( \hat{G}(X^*_{\tau_{n,k}}, \tau_{n,k}) \bigg(\int_t^{\tau_{n,k}} \e^{-\rho s} \drm D^*_s -  \e^{-\rho t} \Delta D^*_t \cdot  \id_{\{(x, t) \in \Pay^o\}} \bigg) \bigg) \\
&\quad - \e^{-2\rho t} (\Delta D^*_t)^2 \, \id_{\{(x, t) \in \Pay^o\}},
\end{align*}
in which we have used the fact that $\hat G$ and $\hat H$ are at least $\mC^1$ with respect to $x$, and $0 < X_s^* \le n$ for all $s \in [t, \tau_{n, k}]$ to deduce that all stochastic integrals involved above have zero expectations.

The growth conditions on $G$ and $H$ imply that as $n \to \infty$ and $k \to \infty$, $\E_{x, t} \big( \hat{H}(X^*_{\tau_{n,k}}, \tau_{n,k}) \big) \to 0$ and $\E_{x, t} \big( \hat{G}(X^*_{\tau_{n,k}}, \tau_{n,k}) \big) \to 0$; also, $\tau_{n,k} \to \tau^*$, and the monotone convergence theorem applies. Therefore, upon sending $n \to \infty$ and $k \to \infty$, we obtain 
\begin{align*} 
\hat{H}(x, t) &= 2  \E_{x, t} \bigg( \int_t^{\tau^*} \e^{-\rho s} \left(  \int_s^{\tau^*} \e^{-\rho u} \drm D^*_u \right) \drm D^*_s \bigg) - \e^{-2\rho t} (\Del D^*_t)^2 \cdot  \id_{\{(x, t) \in \Pay^o\}} \\
&=  2  \E_{x, t} \bigg( \int_{t^+}^{\tau^*} \e^{-\rho s} \left(  \int_s^{\tau^*} \e^{-\rho u} \drm D^*_u \right) \drm D^*_s \bigg) + 2 \e^{-\rho t} \Del D^*_t  \E_{x, t} \bigg( \int_t^{\tau^*} \e^{-\rho u} \drm D^*_u \bigg) \cdot  \id_{\{(x, t) \in \Pay^o\}} \\
&\quad - \e^{-2\rho t} (\Del D^*_t)^2 \cdot  \id_{\{(x, t) \in \Pay^o\}} \\
&= 2 \E_{x, t} \bigg( \int_{t}^{\tau^*} \e^{-\rho s} \left(  \int_s^{\tau^*} \e^{-\rho u} \drm D^{*,c}_u \right) \drm D^{*,c}_s \bigg) \\
&\quad + 
2 \e^{-\rho t} \Del D^*_t \, \E_{x, t} \bigg( \int_t^{\tau^*} \e^{-\rho u} \drm D^{*,c}_u + \e^{-\rho t} \Delta D^*_t \cdot \id_{\{(x, t) \in \Pay^o\}} \bigg)  \cdot  \id_{\{(x, t) \in \Pay^o\}} \\
&\quad - \e^{-2\rho t} (\Del D^*_t)^2 \cdot  \id_{\{(x, t) \in \Pay^o\}} \\
&= \E_{x, t} \bigg( \bigg(\int_t^{\tau^*} \, \e^{- \rho s} \, \drm D^{*,c}_s +  \e^{-\rho t} \Del D^*_t \cdot  \id_{\{(x, t) \in \Pay^o\}} \bigg)^2 \bigg) \\
&= \E_{x, t} \bigg( \bigg(\int_t^{\tau^*} \, \e^{- \rho s} \, \drm D^{*}_s \bigg)^2 \bigg) = \E_{x,t} \big( \left(Y_t^*\right)^2 \big),
\end{align*}
thereby, proving the result in \eqref{eq:GH}.

\medskip
\noi 
\textbf{Step 3.}  We show that if $\Vt$ solves \eqref{eq:HJB_V} with $\Vt(0,t) = 0$, then $\Vt(x, t) = J(x, t; D^*)$.

First, we consider  $(x, t) \in \NT$; in this case, the first term in \eqref{eq:HJB_V} equals 0. Using this identity, along with \eqref{eq:HJB_G} and \eqref{eq:HJB_V}, leads to 
\begin{align*}
	\mM  \Big( \Vt(x, t) - \dfrac{\gam}{2} \, G^2(x, t) \Big) &= - \gam G(x, t) \mM  G(x, t) + \rho G(x, t) - \gam \rho \big(H(x, t) - G^2(x, t) \big) \\
	&= \rho G(x, t) - \gam \rho H(x, t) = \mM  \Big( G(x, t) - \dfrac{\gam}{2} \, H(x, t) \Big).
\end{align*} 
We, then, apply similar arguments as in Step 1 to $\Vt(x, t) - \frac{\gam}{2} \, G^2(x, t)$ and $G(x, t) - \frac{\gam}{2} \, H(x, t)$ and use the above equality, the transversality condition, and $\Vt (0, t) -  \frac{\gam}{2} \, G^2(0, t) = 0 = G(0, t) - \frac{\gam}{2} H(0, t)$ to conclude that  $\Vt(x, t) = G(x, t) - \frac{\gam}{2} ( H(x, t) - G^2(x, t) ) = \E_{x,t}(Y_t^*) - \frac{\gam}{2} \V_{x, t} (Y_t^*) = J(x, t; D^*)$.

Next, we consider $(x, t) \in \Pay$; in this case, $\partial_x \Vt(x, t) = \partial_x G(x, t) = 1$ and $\partial_x H(x, t) = 2 G(x, t)$. As such, 
\[
\partial_x  \Big( \Vt(x, t) - \dfrac{\gam}{2} \, G^2(x, t) \Big) = 1 - \gam G(x, t) = \partial_x  \Big( G(x, t) - \dfrac{\gam}{2} \, H(x, t) \Big),
\]
which, along with the boundary condition, confirms $\Vt(x, t) = J(x, t; D^*)$ for this case. 

Therefore, for all $(x, t) \in \R_+^2$, $\Vt(x, t) = J(x, t; D^*)$ holds as desired. 

\medskip
\noi 
\textbf{Step 4.} We show that if $D^* \in \mA_t$ solves \eqref{eq:Sk}, then $D^*$ is an equilibrium dividend strategy.

To that end, define the perturbed strategy $D^\eps$ as in \eqref{eq:pD_eps}, and we want to prove that the limit in \eqref{eq:pJ_eqm} holds.  Recall that $d \in [0, x]$ is the lump-sum payment at $t$ under $D^\eps$. First, assume $d = x$, and ruin occurs immediately at $t$ under $D^\eps$, resulting in $J(x, t; D^\eps) = x$. Because $V(0, t) = 0$ and $\partial_x V(x, t) \ge 1$, it follows that $V(x, t) \ge x = J(x, t; D^\eps)$ for all $(x, t)$, and, thus, the limit in \eqref{eq:pJ_eqm} holds when $d = x$. 

Given the above analysis, we assume $d < x$ in the remainder of this step; we also  write $X := X^{D^\eps}$, $Y_t := Y_t^{D^\eps}$, and $\tau := \tau^{D^\eps}$ for notational simplicity in the proof. 
By definition, $J(x, t; D^\eps) = \E_{x, t} (Y_t) - \frac{\gam}{2} \E_{x,t} \big( Y_t^2\big) + \frac{\gam}{2} \big( \E_{x,t}  (Y_t) \big)^2$; in what follows, we analyze each of the three terms in $J(x, t; D^\eps)$ by expanding them to order $o(\eps)$. 
To start, we recall an important result on finite-time ruin probabilities (see, for instance, Appendix in Grandell \cite{G1991})
\begin{align}
	\label{eq:ruin_prob}
	\P_{x, t} \big( \tau > t+ \eps \big) \sim 1 - \frac{b \sqrt{\eps}}{x} \, \exp \left( - \frac{1}{2 \eps} \left(\frac{x}{b}\right)^2 \right) = 1 + o(\eps);
\end{align}
with this result, we can omit $\id_{\{\tau > t + \eps\}}$ in the following derivations. 
For convenience, in the derivation below, we introduce 
\begin{align}
		G_t &:= G(x - d, t),   & H_t &:= H(x - d, t); \\
	\mathbf{I}_\rho &:= \int_{t^+}^{t + \eps} \e^{-\rho (s - t)} \drm \delta(s) , 
	& \mathbf{I}_\phi &:= \E_{x, t} \left(\int_{t^+}^{t + \eps} \partial_x \phi(X_s, s) \drm \delta(s)\right), \; \phi \in \mC^{1,1}(\R_+^2). \label{eq:I}
\end{align}
Recall from Definition \ref{def:sing_timeconsist} that  $\delta$ is a non-decreasing,  continuous function over $[t, t+\eps)$, satisfying $\delta(t+\eps) - \delta(t) = O(\eps)$ as $\eps \to 0$. As such, for $\eps$ small enough, $\mathbf{I}_\rho$ can be approximated by
\begin{align*}
	\mathbf{I}_\rho = \int_{t^+}^{t + \eps} 1 \,  \drm \delta(s) + o(\eps) = \mathbf{I}_\phi + o(\eps), \text{ with } \phi(x, s) \equiv x.
\end{align*}
Now using the fact that $\phi$ is at least $\mC^1$ for any of $\phi = x, G$, or $H$, the following estimates hold:
\begin{align}
	\label{eq:order}
	\mathbf{I}_\phi = \mathrm{C}  \eps + o(\eps)
	\quad \text{ and } \quad 
	\mathbf{I}_{\phi} \mathbf{I}_{\phi'} = \mathrm{C}' \eps^{2} + o(\eps^{2}),
	\quad \phi, \phi' = x, G, H,
\end{align}
for some positive constants $\mathrm{C}$ and $\mathrm{C}'$ that might depend on $\phi$ and $\phi'$. Note that the latter result allows us to safely drop terms involving $\mathbf{I}_{\phi} \mathbf{I}_{\phi'}$ if we truncate at the order $o(\eps)$.

First, we analyze $\E_{x, t} (Y_t)$ as follows:
\begin{align*}  
\E_{x, t} (Y_t) &= \E_{x, t}\bigg(\int_t^{\tau} \e^{-\rho (s - t)} \drm D^\eps_s \bigg) 
= \E_{x, t}\bigg(d + \int_{t^+}^{t + \eps} \e^{-\rho (s - t)}  \drm \delta(s) + \int_{t+\eps}^\tau \e^{-\rho (s - t)}  \drm D^*_s \bigg) \notag \\
&= d + \mathbf{I}_\rho + \E_{x, t}\bigg(\e^{-\rho \eps} \id_{\{ \tau > t + \eps\}}  \, \E_{X_{t+\eps}, t+\eps}\bigg(\int_{t+\eps}^{\tau^*} \e^{-\rho (s - (t + \eps))}  \drm D^*_s \bigg) \bigg) + o(\eps)\notag \\
&= d + \mathbf{I}_\rho + (1 -\rho \eps) \, \E_{x, t} \big(G(X_{t+\eps}, t+\eps) \big) + o(\eps) \notag \\
&= d + \mathbf{I}_\rho + (1 -\rho \eps) \,  \, \E_{x, t} \bigg( G(x, t) +   \int_{t^+}^{t+\eps} \mM G(X_s, s) \, \drm s  + \int_{t^+}^{t+\eps}  b \, \partial_x G(X_s, s) \, \drm B_s
\\
& \quad   - \int_{t^+}^{t+\eps} \, \partial_x G(X_s, s) \, \drm \delta(s) + G(x - d, t) - G(x, t) \bigg) + o (\eps)
\\
&= d + G_t + \eps\big( \mM G_t - \rho G_t \big) 
+ \mathbf{I}_\rho - \mathbf{I}_G + o(\eps) .
\end{align*}
Next, we consider $\E_{x,t} \big( (Y_t)^2\big)$; by using the It\^o's expansion for $G(X_{t+\eps}, t+\eps)$ and $H(X_{t+\eps}, t+\eps)$ as above, along with \eqref{eq:GH} and \eqref{eq:ruin_prob}, we obtain
\begin{align*}
	\E_{x,t} \big( Y_t^2\big) &= \E_{x, t} \bigg[ \Big(d + \mathbf{I}_\rho +  \int_{t+\eps}^\tau \e^{-\rho (s - t)}  \drm D^*_s \Big)^2 \bigg] \\
	&= d^2 + \mathbf{I}_\rho^2 + \e^{-2 \rho \eps} \, \E_{x, t} \big(H(X_{t+\eps}, t + \eps)\big) 
	+ 2 d \mathbf{I}_\rho \\
	&\quad + 2 (d + \mathbf{I}_\rho)  \e^{- \rho \eps} \, \E_{x, t} \big(G(X_{t+\eps}, t + \eps)\big) + o(\eps) \\
	&= d^2 + \mathbf{I}_\rho^2 + (1 - 2 \rho \eps) \big(H_t + \eps \, \mM H_t - \mathbf{I}_H \big) 	+ 2 d \mathbf{I}_\rho \\ 
	&\quad + 2 (d + \mathbf{I}_\rho) (1 - \rho \eps) \big(G_t + \eps \, \mM G_t - \mathbf{I}_G \big) + o(\eps) \\
	&= d^2 + \mathbf{I}_\rho^2 + H_t + \eps \left(\mM H_t - 2 \rho H_t \right)  - \mathbf{I}_H + 2 d \mathbf{I}_\rho \\
	&\quad + 2 d \big( G_t + \eps \left(\mM G_t - \rho G_t \right) \big)
	+ 2 \mathbf{I}_\rho G_t - 2 d \mathbf{I}_G - 2 \mathbf{I}_\rho \mathbf{I}_G + o(\eps).
\end{align*}
We proceed to analyze the third term $\big( \E_{x,t}  (Y_t) \big)^2$. By using the results about $\E_{x,t}  (Y_t)$, we get 
\begin{align*}
	\big( \E_{x,t}  (Y_t) \big)^2 &= \big( d + G_t + \eps\big( \mM G_t - \rho G_t \big) 
	+ \mathbf{I}_\rho - \mathbf{I}_G + o(\eps) \big)^2 \\
	&= (d + G_t)^2 + 2 \eps (d + G_t) (\mM G_t - \rho G_t) + 2 (d + G_t) (\mathbf{I}_\rho - \mathbf{I}_G) + (\mathbf{I}_\rho - \mathbf{I}_G)^2 + o(\eps).
\end{align*}

By combining the analysis of the three terms above and using the approximation $\mathbf{I}_\rho = \mathbf{I}_x + o(\eps)$, we obtain 
\begin{align*}
	J(x, t; D^\eps) &= d + G_t + \eps\big( \mM G_t - \rho G_t \big) 
	 \notag \\
	& \quad - \frac{\gam}{2} \big(d^2 + \big( H_t + \eps (\mM H_t - 2 \rho H_t ) \big) + 2 d \big( G_t + \eps \left(\mM G_t - \rho G_t \right) \big)\big) \\
	& \quad + \frac{\gam}{2} \left( (d + G_t)^2 + 2 \eps (d + G_t) (\mM G_t - \rho G_t) \right)
	\\
	& \quad - \frac{\gam}{2} \left( \mathbf{I}_\rho^2 - \mathbf{I}_H + 2 d \mathbf{I}_\rho + 2 \mathbf{I}_\rho G_t - 2 d \mathbf{I}_G - 2 \mathbf{I}_\rho \mathbf{I}_G  \right) \\
	& \quad  + \frac{\gam}{2}  \left( 2 (d + G_t) (\mathbf{I}_\rho - \mathbf{I}_G) + (\mathbf{I}_\rho - \mathbf{I}_G)^2 \right) + \mathbf{I}_\rho - \mathbf{I}_G + o(\eps) 
	\\
	& = d + \left( G_t - \frac{\gam}{2} \left(H_t - G_t^2 \right) \right) 
	+ \eps (1 + \gam G_t) (\mM G_t - \rho G_t) - \frac{\gam}{2} \eps (\mM H_t - 2 \rho H_t ) \\
	&\quad + \mathbf{I}_x - \mathbf{I}_G - \gam G_t \mathbf{I}_G +\frac{\gam}{2} \mathbf{I}_H    + o(\eps).
\end{align*}
Note from the proof of Step 3 that $G - \frac{\gam}{2} \left(H - G^2 \right) = \Vt$. Using the definition of $\mathbf{I}_\phi$, with $\phi = x, G, H$, we have 
\begin{align}
	\label{eq:int_est}
	\mathbf{I}_x - \mathbf{I}_G - \gam G_t \mathbf{I}_G +\frac{\gam}{2} \mathbf{I}_H &=  \int_{t^+}^{t + \eps} \left(1 - \partial_x G - \gam G \partial_x G + \frac{\gam}{2} \partial_x H \right)\Big|_{(X_s, s)} \, \drm \delta(s) + o(\eps) \notag \\
	&= \int_{t^+}^{t + \eps} \left(1 - \partial_x \Vt(X_s, s) \right) \drm \delta(s) + o(\eps) \le o(\eps),
\end{align}
in which the last inequality follows from $ 1 - \partial_x \Vt \le 0$ in \eqref{eq:HJB_V}.
Therefore, by using the above estimate and the identity of $\Vt$, we further reduce $J(x, t; D^\eps)$ to 
\begin{align}
	J(x, t; D^\eps) &= d + \Vt_t + \eps \left( \mM \Vt_t - \frac{\gam}{2}   \mM G^2_t + \gam G_t  \mM G_t - \rho G_t + \gam \rho \big(H_t - G^2_t \big) \right) + o(\eps) \\
	& \le d + \Vt_t + o(\eps) \\
	&= \Vt(x, t) + \left(d - \int_{x -d}^x \, \partial_x \Vt(z, t) \, \drm z \right) + o(\eps) \\
	& \le \Vt(x, t) + o(\eps) = J(x, t; D^*)  + o(\eps),
\end{align}
in which $\Vt_t := \Vt(x - d, t)$, and the two inequalities follow from \eqref{eq:HJB_V}. 

Finally, we conclude that the desired limit result in \eqref{eq:pJ_eqm} holds. 
\end{proof}

\begin{remark}
	\label{rem:veri}
	Because both this paper and Cao et al.\  {\rm\cite{CLYZ2025}} seek equilibrium strategies for MV dividend problems from a game-theoretic perspective, their verification lemmas and proofs share certain similarities. However, there are also major differences, which we now describe. In Theorem {\rm \ref{thm:verif_sing}}, the feasible region is partitioned into the pay region and no-transaction region, and all related functions $(V$, $G$, and $H)$ are characterized separately in these two regions; the value function $V$ satisfies an HJB-variational inequality equation in \eqref{eq:HJB_V}, which appears because of the singular-control framework. By comparison, each corresponding function in Theorem $2.3$ of Cao et al.\ {\rm \cite{CLYZ2025}} is characterized by one second-order PDE over the entire region, and the value function $V$ therein satisfies a standard HJB equation. Regarding the equilibrium strategies, $D^*$ in this paper is obtained as a solution of the associated Skorokhold reflection problem in \eqref{eq:Sk}, but $D^*$ in  Cao et al.\ {\rm \cite{CLYZ2025}} is the maximizer of the HJB equation of $V$. The difference in their proofs lies on the technical side and mainly originates from their differences concerning the definition of perturbed strategies $($see Remark {\rm \ref{rem:eq_def}}$)$. In particular, it  takes a delicate and involved analysis to study the performance of the perturbed strategy $D^\eps$ here in Step $4$, which eventually yields the desired first-order inequality in \eqref{eq:pJ_eqm}. \qed
\end{remark}

\section{Equilibrium dividend strategies}\label{sec:eqm}

In this section, we apply the verification theorem (Theorem \ref{thm:verif_sing}) to derive the equilibrium dividend strategy $D^*$ in (semi)closed form for large $\gamma$ (risk aversion) and small $\gamma$. 

To begin, we review the special case of $\gam = 0$; note that the objective $J$ in \eqref{eq:J} becomes $\E_{x, t}(Y_t)$, and the corresponding optimal dividend problem is time-consistent and has been solved in the literature. For instance, Theorem 2.2 and Lemma 2.3 in Taksar \cite{T2000} show that the optimal strategy is a barrier strategy with a strictly positive barrier $\xt$ (because $a > 0$), and the value function is concave and obtained explicitly in a two-piece form separated by the barrier $\xt$. 
We hypothesize that for small positive $\gam$, a similar result holds. However, for $\gam$ large enough, the penalty on the variation of dividend payments should ``force'' the manager to pay the entire surplus and declare bankruptcy (yielding a zero variance).  We formally verify this latter hypothesis in the next theorem.

\begin{theorem}\label{thm:eqm_pV_eq_x}
If the following condition holds
\begin{equation}\label{eq:gam_large}
	\gam \ge \dfrac{2a}{b^2},
\end{equation}
then an equilibrium dividend strategy is to pay out all of surplus as dividends immediately $($that is, $D_t^* = X_{t^-} = x$ and $\tau^* = t)$, and we have $V(x, t) = G(x, t) \equiv x$ and $H(x, t) \equiv x^2$ for all $x \ge 0$. 
\end{theorem}

\begin{proof}
Suppose inequality \eqref{eq:gam_large} holds; we consider the strategy of paying all dividends immediately (and thereby ruining immediately). This strategy is clearly admissible by Definition \ref{def:adm}, and it implies that the pay region is $\Pay = \R_+^2$ and $Y_t = x$. As such, it follows from \eqref{eq:GH} that $G(x, t) \equiv x$ and $H(x, t) \equiv x^2$, and they satisfy the related HJB equations in \eqref{eq:HJB_G} and \eqref{eq:HJB_H}. Given $G$ and $H$, we obtain the (candidate) value function by $V(x, t) = G(x, t) - \frac{\gam}{2} (H(x, t) - G^2(x,t)) \equiv x$, which implies that $1 - \partial_x V \le 0$ holds with equality, and the boundary condition in \eqref{eq:boundary} is satisfied. It remains to show that the first variational inequality in \eqref{eq:HJB_V} is true for all $x \in \R_+$. To that end, we compute
\begin{align}
\mM V - \dfrac{\gam}{2} \, \mM G^2 + \gam G \mM G - \rho G + \gam \rho \big(H - G^2 \big) 
= a - \dfrac{\gam}{2} \, b^2 - \rho x,
\end{align} 
which is non-positive for all $x \ge 0$ when \eqref{eq:gam_large} holds.  Thus, $V$, $G$, and $H$ satisfy the conditions of Theorem \ref{thm:verif_sing}, and paying out all of surplus as dividends is an equilibrium strategy.
\end{proof}

We next prove a non-trivial result for small risk aversion $\gam$ and confirm the earlier hypothesis that a barrier strategy is an equilibrium dividend strategy. For convenience, define 
\begin{align}
	\label{eq:r_12}
	r_1 &= \dfrac{1}{b^2} \left[  - a + \sqrt{a^2 + 2 \rho b^2} \right] > 0, & 
	r_2 &= \dfrac{1}{b^2} \left[  - a - \sqrt{a^2 + 2 \rho b^2} \right] < 0, 
	\\
	r_3 &= \dfrac{1}{b^2} \left[  - a + \sqrt{a^2 + 4 \rho b^2} \right] > 0, & 
	r_4 &= \dfrac{1}{b^2} \left[  - a - \sqrt{a^2 + 4 \rho b^2} \right] < 0. \label{eq:r_34}
\end{align}

\begin{theorem}\label{thm:onebarrier}
There exists an $\eps \in (0, \frac{2a}{b^2})$ such that for all $\gam \in (0, \eps)$,  equation
\begin{align}\label{eq:xt}
	0 &= \dfrac{r_1^2 \e^{r_1 x} - r_2^2 \e^{r_2 x}}{r_1 \e^{r_1 x} - r_2 \e^{r_2 x}} + \gam \left\{ 1 + \dfrac{\e^{r_1 x} - \e^{r_2 x}}{r_1 \e^{r_1 x} - r_2 \e^{r_2 x}} \left( \dfrac{r_1^2 \e^{r_1 x} - r_2^2 \e^{r_2 x}}{r_1 \e^{r_1 x} - r_2 \e^{r_2 x}} - \dfrac{r_3^2 \e^{r_3 x} - r_4^2 \e^{r_4 x}}{r_3 \e^{r_3 x} - r_4 \e^{r_4 x}} \right)  \right\} =: f(x, \gamma) \qquad 
\end{align}
admits a unique positive solution, denoted by $\tilde{x}$.  If  $(1)$ $\gam < \eps$, and $(2)$ $V$ in \eqref{eq:Vx} is strictly concave over $[0, \xt)$, then a 
barrier strategy, with constant barrier $\xt$, is an equilibrium strategy $D^*$, with $\Pay = [\xt, \infty) \times \R_+$ and $\NT = [0, \xt) \times \R_+$. Moreover,
\begin{align}\label{eq:pG3}
G(x) = \E_{x}(Y_0^*) = 
\begin{cases}
C_1 \big(\e^{r_1 x} - \e^{r_2 x} \big),  &\quad x < \xt, \\
C_1 \big(\e^{r_1 \xt} - \e^{r_2 \xt} \big) + (x - \xt),  &\quad x \ge \xt,
\end{cases}
\end{align}
and
\begin{align}\label{eq:pH3}
H(x) = \E_{x} \big[ (Y^*_0)^2 \big] = 
\begin{cases}
C_3 \big(\e^{r_3 x} - \e^{r_4 x} \big),  &\quad x < \xt, \\
C_3 \big(\e^{r_3 \xt} - \e^{r_4 \xt} \big) + 2 C_1\big(\e^{r_1 \xt} - \e^{r_2 \xt} \big)(x - \xt) + (x - \xt)^2,  &\quad x \ge \xt,  \hspace{2em}
\end{cases}
\end{align}
in which 
\begin{align}
	\label{eq:C_13}
	C_1 =  \dfrac{1}{r_1 \e^{r_1 \xt} - r_2 \e^{r_2 \xt}} > 0 \quad \text{ and } \quad 
	C_3 = \dfrac{2\big(\e^{r_1 \xt} - \e^{r_2 \xt} \big)}{ \left( r_1 \e^{r_1 \xt} - r_2 \e^{r_2 \xt} \right) \left( r_3 \e^{r_3 \xt} - r_4 \e^{r_4 \xt} \right) } > 0,
\end{align}
and the corresponding value function equals 
\begin{align}
	\label{eq:Vx}
	V(x) = G(x) - \frac{\gam}{2} \left( H(x) - G^2(x) \right).
\end{align}
\end{theorem}

\begin{proof}
Because the problem is time-homogeneous, we expect the value function $V$ to be time-independent, along with $G$ and $H$. For this reason, we set time equal to $0$ and suppress the time argument in the analysis; also, we write $\phi'$ and $\phi''$ to denote the first and second derivative (with respect to $x$) for $\phi = V, G$, or $H$. 
We hypothesize that a time-independent barrier strategy, with a constant barrier $\xt > 0$, is an equilibrium dividend strategy $D^*$ (in the sense of Definition \ref{def:sing_timeconsist}). Specifically, this strategy dictates the manager of the company to pay $(x - \xt) \id_{x \ge \xt}$ in dividends at time 0 (with initial surplus $X_{0^-} = x \ge 0$) and thereafter pay dividends in order to keep the surplus $X_t^* \in [0, \xt]$ for all $t > 0$. Since the barrier strategy $D^*$ is time-independent, we write the pay region as $\Pay = [\xt, \infty)$ and the no-transaction region as $\NT = [0, \xt)$ associated with $D^*$ in the proof.

With the above hypothesis, we proceed to solve for $G$ and $H$ based on whether $x \in \NT$ or $x \in \Pay$. 
First, assume $x \in \NT$, that is, $x < \xt$. In this case, by \eqref{eq:HJB_G}, $G$ solves the boundary-value problem,
 $- \rho G(x) + a G'(x) + \frac{1}{2} \, b^2 G''(x) = 0$, with $G(0) = 0$, 
whose solution equals the first expression in \eqref{eq:pG3}. Similarly, using \eqref{eq:HJB_H} for $H$, we solve $- 2 \rho H(x) + a H'(x) + \frac{1}{2} \, b^2 H''(x) = 0$, given $H(0) = 0$, 
and obtain the first expression for $H$ in \eqref{eq:pH3}.

Next, assume $x \in \Pay$, that is, $x \ge \xt$.  In this case, our ansatz strategy implies that the company immediately pays a lump-sum dividend of $x - \xt$. By the continuity of $G$ and using $G'(x) = 1$, we arrive at the second expression of $G$ in \eqref{eq:pG3}. Next, \eqref{eq:HJB_H} implies $H'(x) = 2 G(x)$, and using this result leads to the second expression of $H$ in \eqref{eq:pH3}.

The two positive constants $C_1$ and $C_3$ in \eqref{eq:pG3} and \eqref{eq:pH3} are yet to be determined. To determine them, we use the condition that $G, H \in \mC^2(\R_+)$, except possibly at $x = \xt$ where they must be $\mC^1$. This motivates us to impose the ``smooth pasting'' condition: $G'(\xt^-) = G'(\xt^+)$ and $H'(\xt^-) = H'(\xt^+)$,
from which we obtain $C_1$ and $C_3$ as in \eqref{eq:C_13}. 

With $G$ and $H$ obtained in \eqref{eq:pG3} and \eqref{eq:pH3}, respectively, we immediately obtain the candidate value function $V$ by \eqref{eq:Vx}; note that the barrier $\xt > 0$ (appearing in \eqref{eq:pG3}, \eqref{eq:pH3}, and \eqref{eq:Vx}) is unknown from the ansatz. 
To determine $\xt$, we impose the condition 
\begin{align}
	\label{eq:V_cond}
	V''(\xt^-) := \lim_{x \uparrow \xt} V''(x)= 0.
\end{align}
By using \eqref{eq:pG3} and \eqref{eq:pH3}, we verify that $V''(\xt^+) :=\lim_{x \downarrow \xt} V''(x) = 0$ holds automatically. Therefore, with \eqref{eq:V_cond}, we have $V''(\xt) = 0$, and it further implies $V \in \mC^2(\R_+)$ because the continuity of $V$ and $V'$ follows from that of $G$, $G'$, $H$, and $H'$.
To obtain a finer condition for \eqref{eq:V_cond}, we compute: for all $x < \xt$,
\begin{align}
V''(x) &= \dfrac{r_1\e^{r_1x} - r_2\e^{r_2x}}{r_1\e^{r_1\tilde{x}} - r_2\e^{r_2\tilde{x}}}\cdot \Bigg\{ \dfrac{r_1^2\e^{r_1x}  - r_2^2\e^{r_2x}}{r_1\e^{r_1x} - r_2\e^{r_2x}} + \gam \bigg[ \dfrac{r_1\e^{r_1x} - r_2\e^{r_2x}}{ r_1\e^{r_1\xt} - r_2\e^{r_2\xt}} + \dfrac{ (\e^{r_1x} - \e^{r_2x})(r_1^2\e^{r_1x} - r_2^2\e^{r_2x})}{(r_1\e^{r_1\xt} - r_2\e^{r_2\xt})(r_1\e^{r_1x} - r_2\e^{r_2x})} \notag \\
&\qquad\qquad - \dfrac{(\e^{r_1\xt} - \e^{r_2\xt})(r_3^2\e^{r_3x} - r_4^2\e^{r_4x})}{(r_3\e^{r_3\xt} - r_4\e^{r_3\xt})(r_1\e^{r_1x} - r_2\e^{r_2x})} \bigg]\Bigg\} =: \dfrac{r_1\e^{r_1x} - r_2\e^{r_2x}}{r_1\e^{r_1\tilde{x}} - r_2\e^{r_2\tilde{x}}}\cdot g(x, \xt),  \label{eq:Vpp}
\end{align}
which shows that $V''(\xt^-) =0$ in \eqref{eq:V_cond} is equivalent to $f(\xt, \gam) = 0$ in  \eqref{eq:xt}.

To study the solvability of \eqref{eq:xt}, we treat the right side of \eqref{eq:xt} as a function of $x$ and $\gam$ and denote it by $f(x, \gam)$. For every fixed $\gam$ satisfying  $\gam < 2a / b^2$, we have 
\begin{align}
	\label{eq:f_lim}
	f(0, \gam) = \gam - \dfrac{2a}{b^2} < 0 \quad \text{ and } \quad \lim_{x \to \infty} f(x, \gam) = r_1 + \gam \left( 2 - \dfrac{r_3}{r_1} \right) > 0,
\end{align}
in which  the second inequality follows from $2 r_1 > r_3$ by their definitions in \eqref{eq:r_12} and \eqref{eq:r_34}. As such, combining with the fact that $f(\cdot, \gam)$ is continuous over $\R_+$, $f(x, \gam) = 0$ admits at least one positive solution $x_\gam$ for all $\gam \in (0, 2a/b^2)$. To obtain the uniqueness result, we first set $\gam = 0$ and verify that $f(x, 0) = 0$ has a unique positive solution, $x_0$ (by using \eqref{eq:f_lim} and verifying $\partial_x f(x, 0) > 0$).  Moreover, by a tedious calculation, we deduce $\partial_x f(x, \gam)|_{(x_0, 0)} \propto - r_1 r_2 (r_1 - r_2)^2 \, \e^{(r_1 + r_2) x_0} > 0$. Therefore, by the implicit function theorem, there exists a small $\eps \in (0, \frac{2a}{b^2})$ such that \eqref{eq:xt} has a unique positive solution $\xt := x_\gam$ (that is, $f(x_\gam, \gam) = 0$) for all $\gam < \eps$. (Recall that we assumed $\gam < \frac{2a}{b^2}$ to obtain $f(0, \gam)<0$ in \eqref{eq:f_lim}; therefore, we impose an upper bound of $\frac{2a}{b^2}$ on $\eps$.)

By construction, $G$ in \eqref{eq:pG3} and $H$ in \eqref{eq:pH3} satisfy all the conditions of Theorem \ref{thm:verif_sing}, and the candidate barrier strategy is admissible and solves the Skorokhod reflection problem \eqref{eq:Sk} with $X_{0^-} = x$ and $\NT = [0, \xt) \times \R_+$ (see, for instance, Lemma 4.1 in Wang and Zou \cite{WZ2021}). The remaining task is to verify that $V$ satisfies the HJB variational equation in \eqref{eq:HJB_V} and that the partition in \eqref{eq:P} is consistent with $V$.  
Because $1 - V'(x) = 0$ on $\Pay$ by \eqref{eq:HJB_G} and \eqref{eq:HJB_H}, the strict concavity of $V$ in Condition (3) implies that $1 - V'(x) < 0$ for all $ x \in \NT$.

Finally, by applying Theorem \ref{thm:verif_sing}, all the results in Theorem \ref{thm:onebarrier} follow as desired.
\end{proof}

\begin{remark}
	\label{rem:gam}
The first condition in Theorem {\rm \ref{thm:onebarrier}} explicitly requires small $\gam$. 
We claim that the second condition $($that is, $V$ in \eqref{eq:Vx} is strictly concave over $[0, \tilde{x}))$ also requires small $\gam$ $(\gam \le \frac{2a}{b^2}$, to be precise$)$. To prove our claim, we argue by contradiction and choose a $\gam > \frac{2a}{b^2}$; suppose $\gam$ is such that $f(x, \gam) = 0$ has a unique solution $\xt$. 
For such a $\gam$, we have $f(0, \gam) > 0$ and $\lim_{x\to\infty}f(x, \gam) > 0$.  From the continuity of $f(\cdot, \gam)$ and uniqueness of $\xt$, we deduce $f(x, \gam) > 0$ for all $x \neq \xt$. This, along with \eqref{eq:Vpp}, implies that $g(x, \xt) >0$ for all $x \neq \xt$, which in turn yields $V''(x) > 0$ over $[0, \xt)$, contradicting the strictly concavity of $V$. 
Numerical analysis in the next section further suggests that there exists an upper bound $\bar{\gam} < \eps \le \frac{2a}{b^2}$ such that both conditions in Theorem {\rm \ref{thm:onebarrier}} hold. As a consequence, when $\gam > \frac{2a}{b^2}$, Theorem {\rm \ref{thm:eqm_pV_eq_x}} shows that paying the entire surplus immediately is an equilibrium strategy; when $\gam \le \bar{\gam}$, Theorem {\rm \ref{thm:onebarrier}} shows that a barrier strategy with a constant barrier $\xt$ is an equilibrium strategy. However, for intermediate level risk aversion $\gam \in (\bar{\gam}, \frac{2a}{b^2})$, finding equilibrium strategies remains an open question. \qed
\end{remark}

\begin{remark} 
	\label{rem:comp}
In this remark, we first compare our results in Theorems {\rm \ref{thm:eqm_pV_eq_x}} and {\rm \ref{thm:onebarrier}} with those in Cao et al.\ {\rm \cite{CLYZ2025}}. In their paper, the model must satisfy a key inequality $($equation $(3.3))$ first, and, then, when risk aversion is small enough, a barrier strategy is an equilibrium strategy $($see Theorem 3.2 in that paper$)$; if the inequality fails, Theorem $3.3$ therein shows that paying dividends at the maximum rate is an equilibrium strategy, but again for \emph{small} risk aversion. By comparison, paying out all surplus in Theorem {\rm \ref{thm:eqm_pV_eq_x}} is an equilibrium strategy for \emph{large}, not small, risk aversion; moreover, a similar inequality is \emph{not} needed for either theorem here. Recall that Condition $(2)$ in Theorem {\rm \ref{thm:onebarrier}} helps verify $1 - V'(x) < 0$ for all $x \in \NT$, which arises from the variational inequality in \eqref{eq:HJB_V}, but a similar condition is \emph{not} needed in Theorem $3.2$ of Cao et al.\ {\rm \cite{CLYZ2025}} because they adopt the classical control framework and the value function only needs to satisfy an $($extended$)$ HJB equation $($see equation $(2.3)$ in that paper$)$.	
	
Under the same diffusion model as ours in  \eqref{eq:X},  Grandits et al.\ {\rm \cite{GHSZ2007}} investigate the optimal dividend strategy that maximizes the expected exponential utility of total dividends paid up to the ruin time, $\max_{D}\E[ U(\int_t^\tau \e^{-\rho s}\drm D_s)]$, in which $U(x) = (1-\e^{-\gam x})/\gam$. They show that  when $\gam \ge \frac{2a}{b^2}$, the optimal strategy is to pay out the entire surplus immediately, which aligns with our finding in Theorem {\rm \ref{thm:eqm_pV_eq_x}}. When $\gam < \frac{2a}{b^2}$, and assuming the existence of a positive solution $b(t)$ to a certain integral equation $($equation $(25))$, the  barrier strategy with time-dependent barrier $b(t)$ is an optimal strategy, which resembles our result in Theorem {\rm \ref{thm:onebarrier}}. The reason that their barrier is time-dependent is that future dividends are discounted to time $0$ in their time-$t$ value function $($equation $12)$; the same setup is also used in Eisenberg and Kr\"uhner {\rm \cite{EK2023}}.  
Gerber and Shiu {\rm \cite{GS2004}} provide a detailed study on the distribution of $Y_0 = \int_0^\tau \, \e^{- \rho s} \, \drm D_s$ under barrier strategies, but they do not attempt to solve for the optimal barrier. \qed
\end{remark}

\section{Numerical examples}
\label{sec:nume}

When risk aversion is large enough ($\gam \ge \frac{2a}{b^2}$), Theorem \ref{thm:eqm_pV_eq_x} shows that paying all of surplus immediately ($D^*_t = x$) is an equilibrium strategy, and the corresponding value function is $V(x) = x$. However, for small $\gam$, the results in Theorem \ref{thm:onebarrier} are less explicit; thus, the first objective of this section is to offer more insights via a detailed numerical analysis. To that end, we set $a = 1$ (surplus drift), $b = 0.25$ (surplus volatility), and $\rho = 0.2$ (discount rate). When $\gamma = 0$, the unique barrier $\xt$ is given by (see equation (2.25) in Taksar \cite{T2000})
\begin{align*}
	\xt = \frac{b^2}{\sqrt{a^2 + 2 \rho b^2}} \, \ln \frac{\sqrt{a^2 + 2 \rho b^2} + a}{\sqrt{a^2 + 2 \rho b^2} - a} \; (=0.3141),
\end{align*}
and the value function $V$ equals $G$ in \eqref{eq:pG3} with the above $\xt$. In this case, $V$ is strictly concave over the NT region and linear over the Pay region. 

\begin{figure}[h]
	\begin{subfigure}{0.5\textwidth}
		\includegraphics[width = 8cm, height = 5cm, trim = 1cm 0cm 1cm 0.8cm, clip = true]{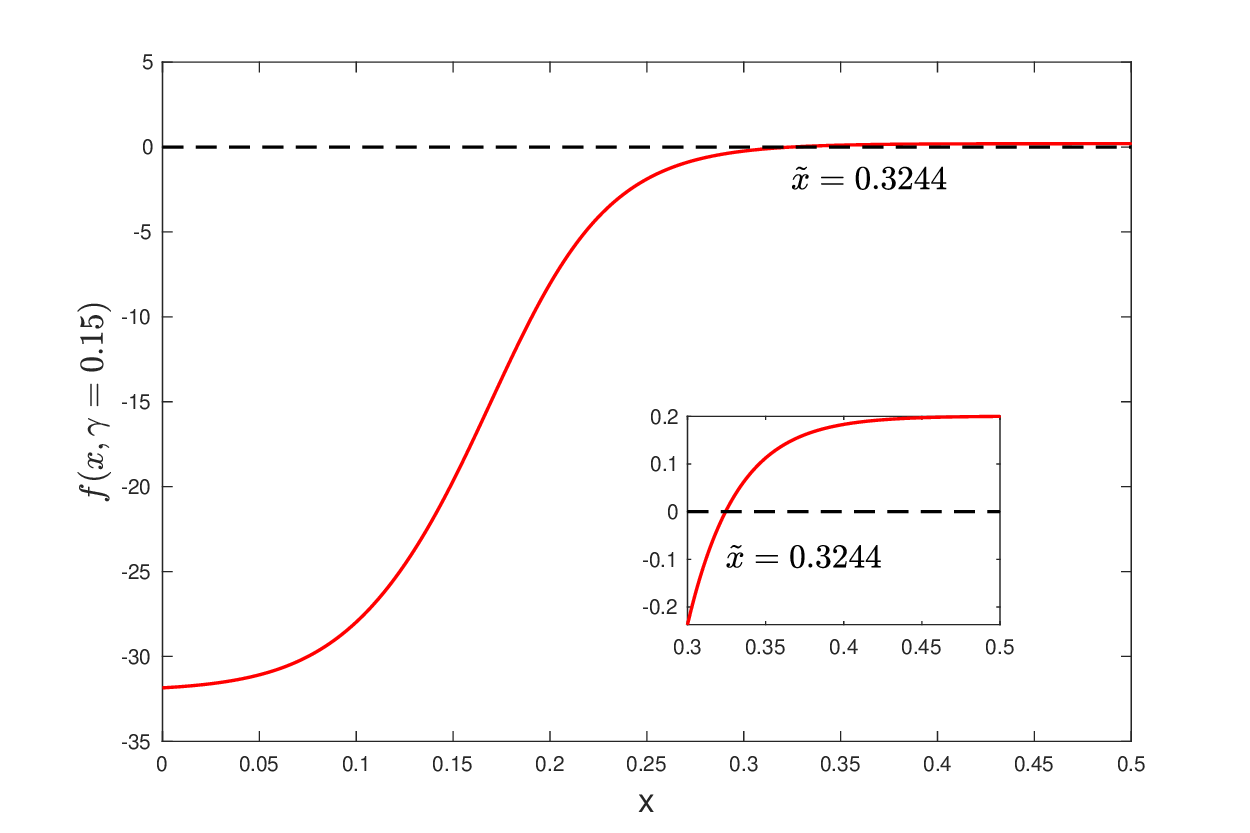}
	\end{subfigure}
	\begin{subfigure}{0.5\textwidth}
		\includegraphics[width = 8cm,  height = 5cm, trim = 1cm 0cm 1cm 0.8cm, clip = true]{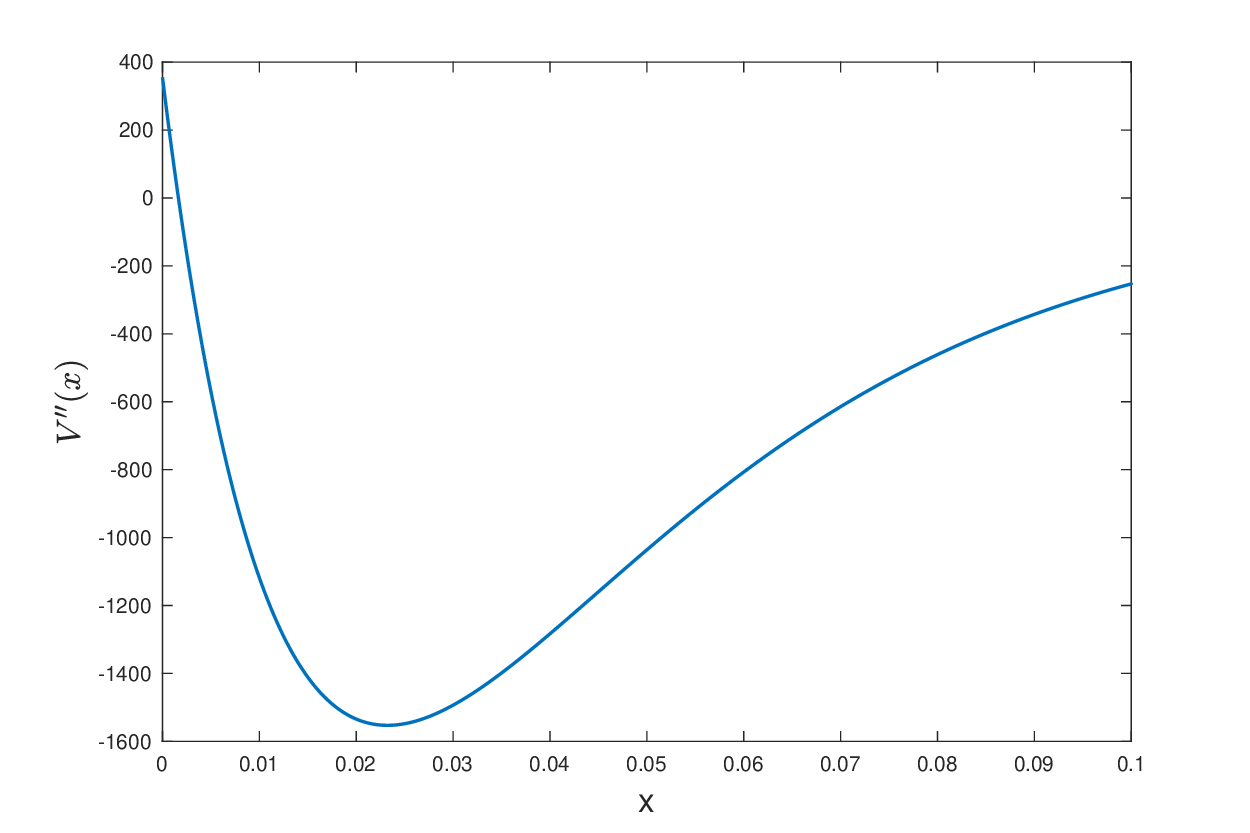}
	\end{subfigure}
\vspace{-4ex}
	\caption{$f(x, 0.15)$ in \eqref{eq:xt} (left) and $V''(x)$ defined via \eqref{eq:Vx} (right) when $\gam = 0.15$}
	\label{fig:gam015}
\end{figure}

\begin{figure}[h]
	\begin{subfigure}{0.5\textwidth}
		\includegraphics[width = 8cm, height = 5cm, trim = 1cm 0cm 1cm 0.8cm, clip = true]{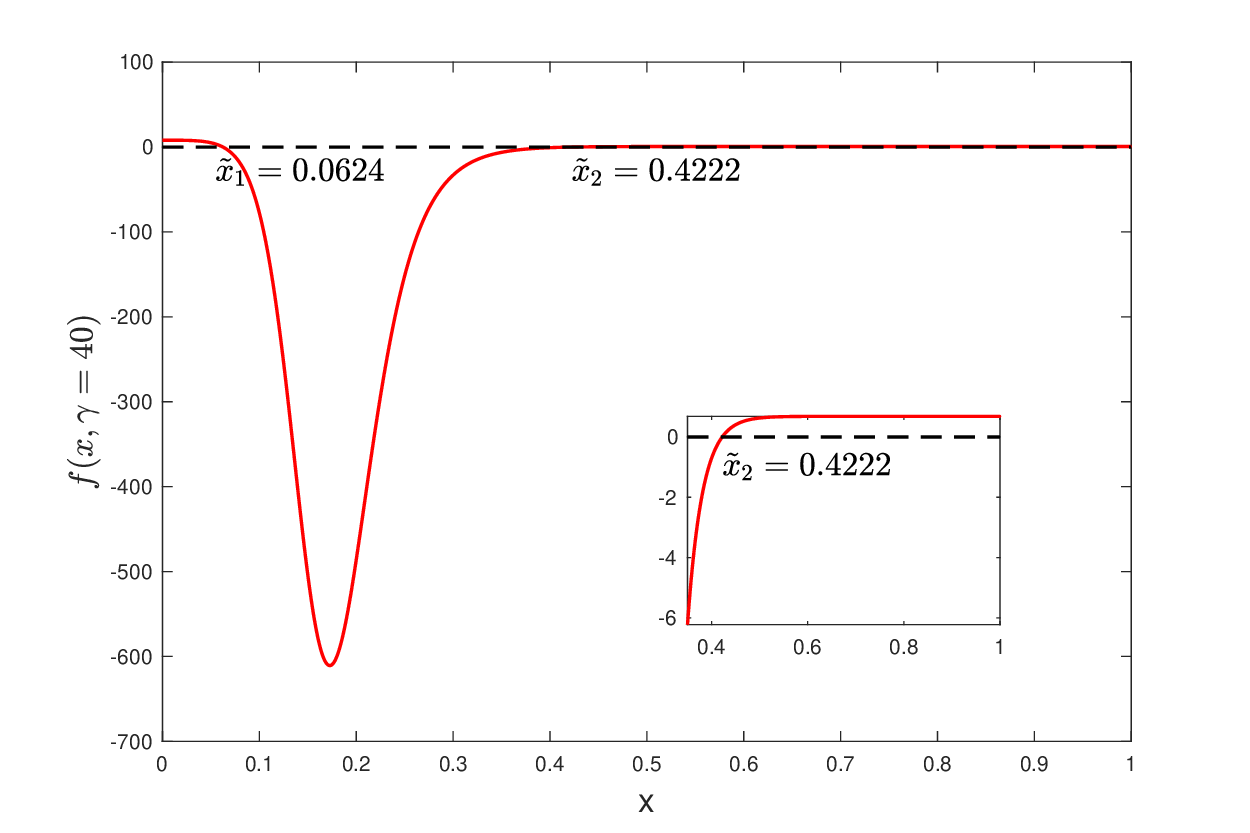}
	\end{subfigure}
	\begin{subfigure}{0.5\textwidth}
		\includegraphics[width = 8cm, height = 5cm, trim = 1cm 0cm 0.5cm 0.8cm, clip = true]{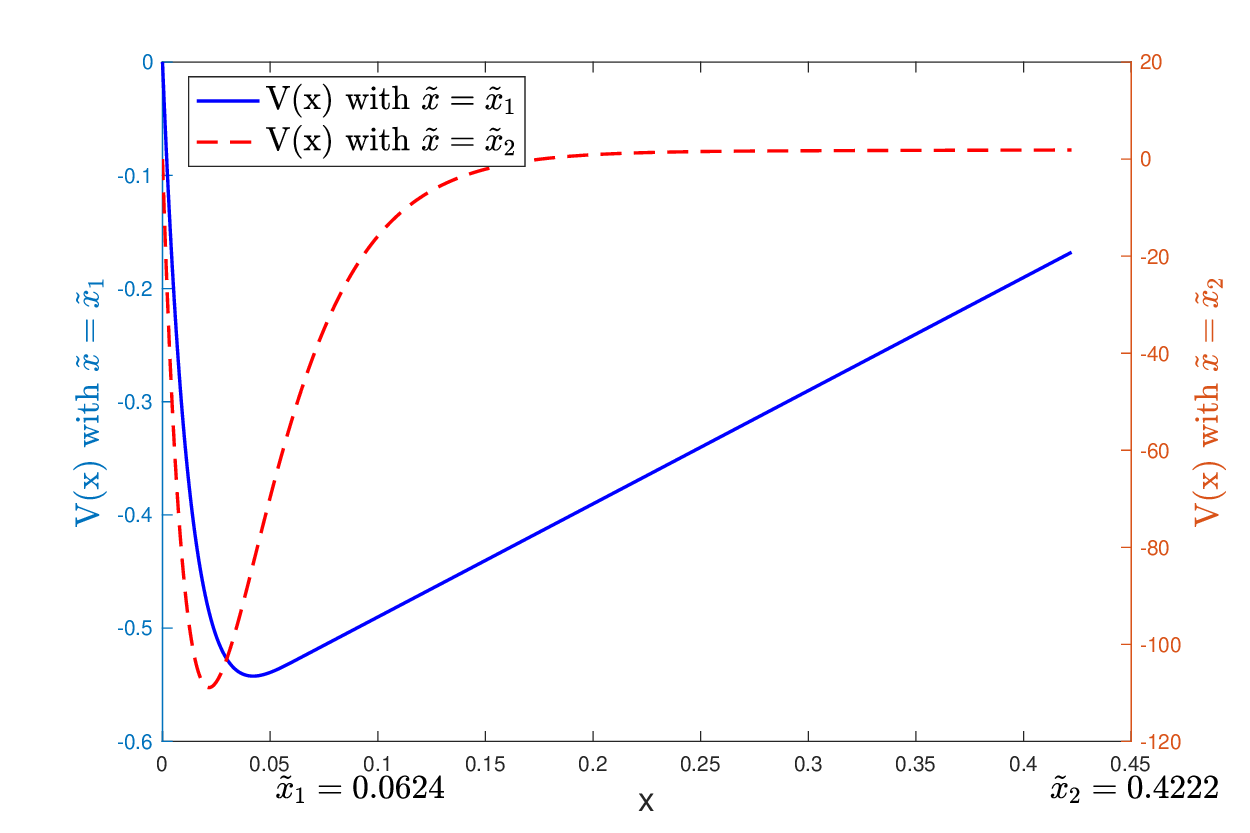}
	\end{subfigure}
\vspace{-4ex}
	\caption{$f(x, 40)$ in \eqref{eq:xt} (left) and the corresponding $V$s in \eqref{eq:Vx} when $\gam = 40$}
	\label{fig:gam40}
\end{figure}

For the given parameter values, our numerical algorithm finds a unique positive root $\xt$ to \eqref{eq:xt} for all $\gam < 32 \, (= \frac{2a}{b^2})$. However, upon substituting $\xt$ into the expression of $V$ in \eqref{eq:Vx}, we find that $V$ is strictly concave over $[0, \xt)$ only for $\gamma \leq 0.1397$. (Strict concavity is required by Condition (2) in Theorem \ref{thm:onebarrier}.)  For instance, if we set  $\gam = 0.15$, then the unique root is $\tilde{x} = 0.3244$, but Figure \ref{fig:gam015} clearly shows that $V''(x) > 0$ when $x$ is near 0. On the other hand, consider $\gam = 40 > \frac{2a}{b^2}$; Figure \ref{fig:gam40} shows that $f(x, 40) = 0$ has two positive roots, $\tilde{x}_1= 0.0624$ and $\tilde{x}_2 = 0.4222$, and plots their corresponding $V$s defined by \eqref{eq:Vx} for $\xt = \xt_1, \, \xt_2$, neither of which is concave over $[0, \xt)$. To echo our earlier comment from Remark \ref{rem:gam}, paying out full surplus is an equilibrium strategy for all $\gam \ge \frac{2a}{b^2} = 32$; a barrier strategy with barrier $\xt$ is an equilibrium strategy for all $\gam \le \bar{\gam} = 0.1397$. But, for intermediate values $\gamma \in (0.1397, 32)$, finding equilibrium strategies remains an open question.

\begin{figure}[h]
	\begin{subfigure}{0.5\textwidth}
		\includegraphics[width = 8cm, height = 5cm, trim = 1cm 0cm 1cm 0.8cm, clip = true]{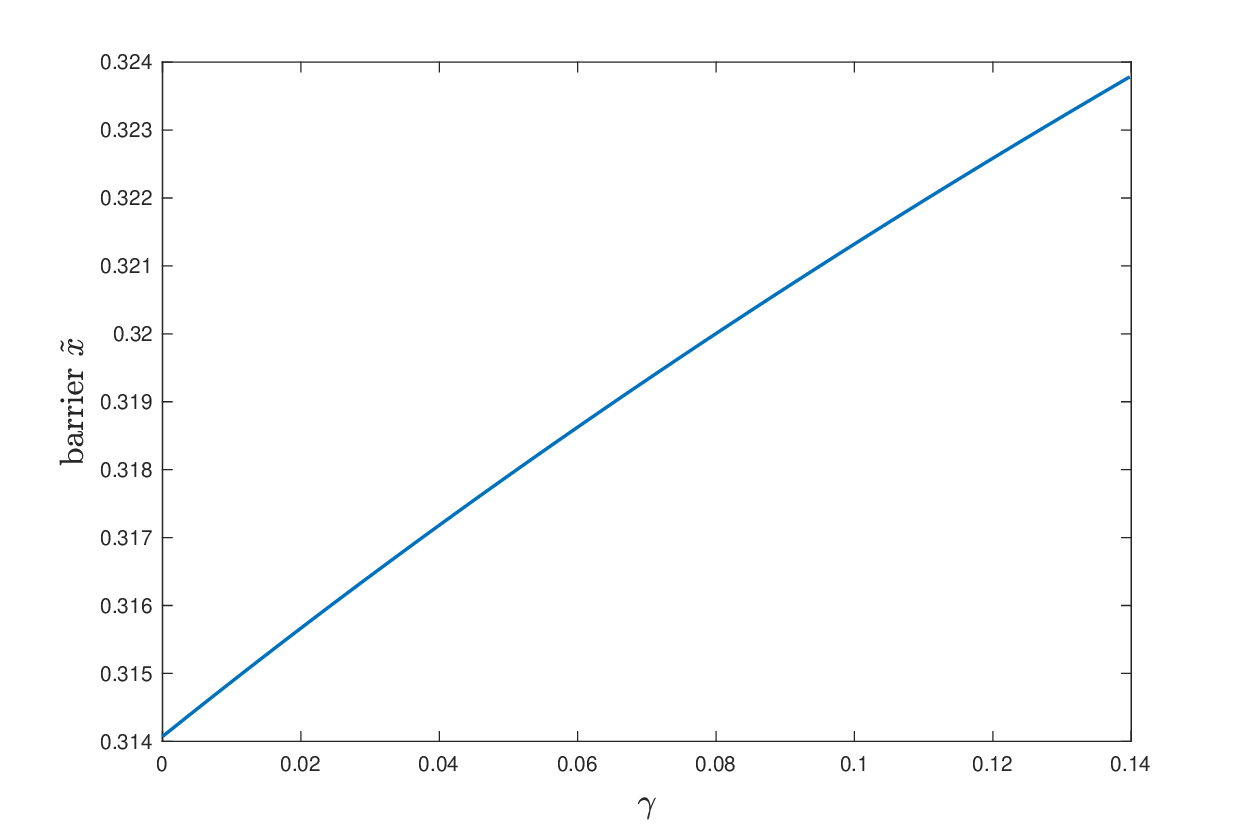}
	\end{subfigure}
    \begin{subfigure}{0.5\textwidth}
    	\includegraphics[width = 8cm, height = 5cm, trim = 1cm 0cm 1cm 0.8cm, clip = true]{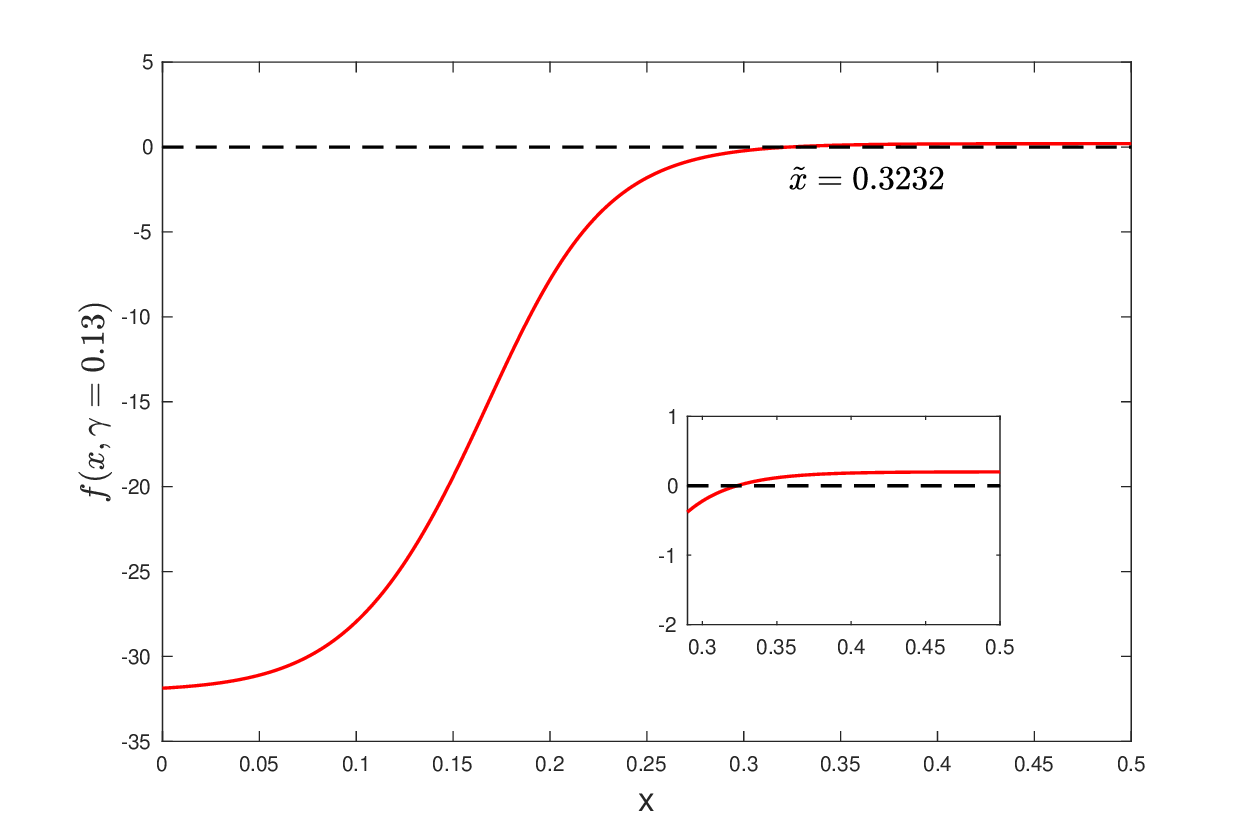}
    \end{subfigure}
\vspace{-4ex} 
\caption{The barrier $\xt$ as a function of $\gamma$ (left) and $f(x, 0.13)$ in \eqref{eq:xt} when $\gam = 0.13$ (right)}
\label{fig1}
\end{figure}

Given the above results, we now consider small $\gam \le 0.1397$ and plot the unique positive root $\xt := \xt(\gam)$ as a function of $\gam$ over this range. The left panel of  Figure \ref{fig1} shows that $\xt$ is an increasing function of $\gam$, and we explain this finding as follows: when $\gamma$ increases, the penalty on dividend variability increases, but since $\gam$ remains small, barrier strategies are still equilibrium strategies; the combined effect, then, drives the manager to set a higher barrier for paying dividends to reduce the variance. The right panel of Figure \ref{fig1} verifies that $f(\cdot, \gam = 0.13) = 0$ has a unique root at $\xt = 0.3232$ when $\gam = 0.13$ (note that $f(\cdot, \gamma)$ is strictly increasing).

\begin{figure}[h]
	\begin{subfigure}{0.5\textwidth}
		\includegraphics[width = 8cm, height = 5cm, trim = 1cm 0cm 1cm 0.8cm, clip = true]{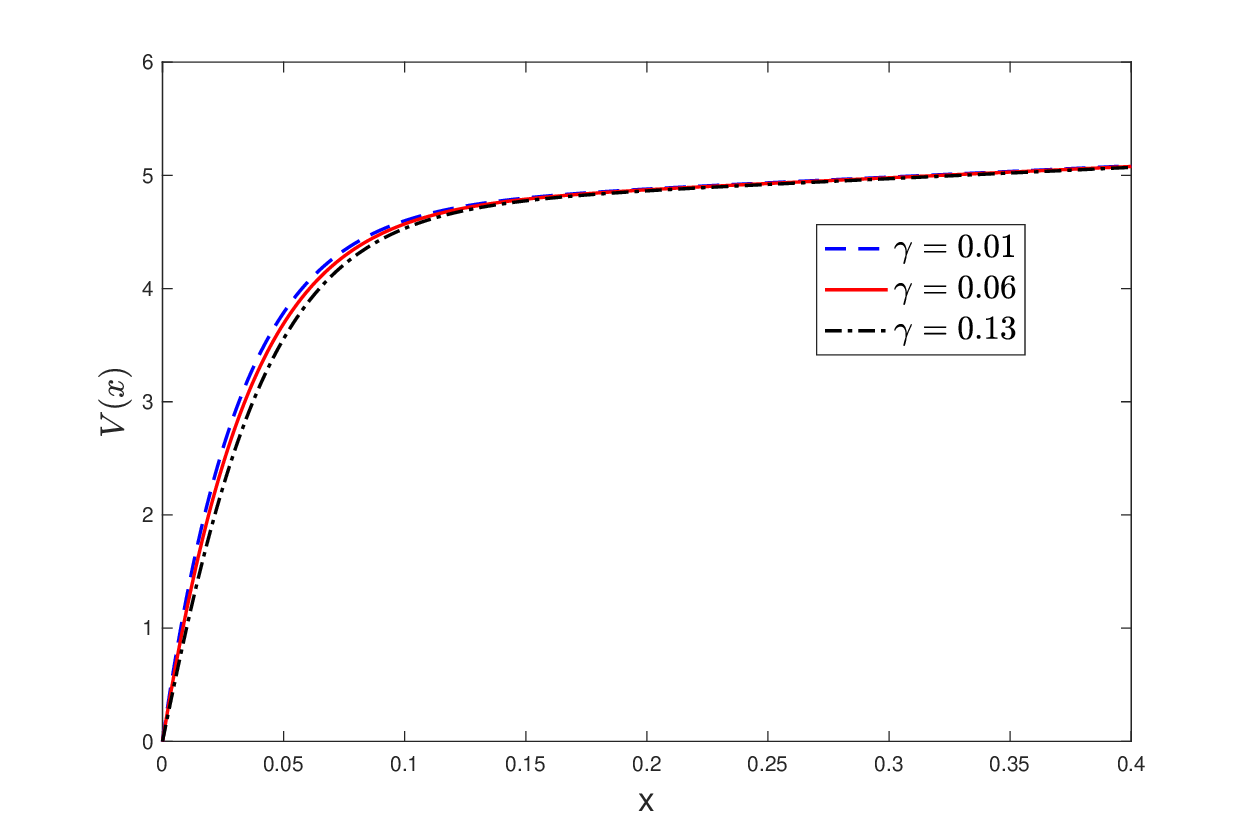}
	\end{subfigure}
	\begin{subfigure}{0.5\textwidth}
		\includegraphics[width = 8cm, height = 5cm, trim = 1cm 0cm 1cm 0.8cm, clip = true]{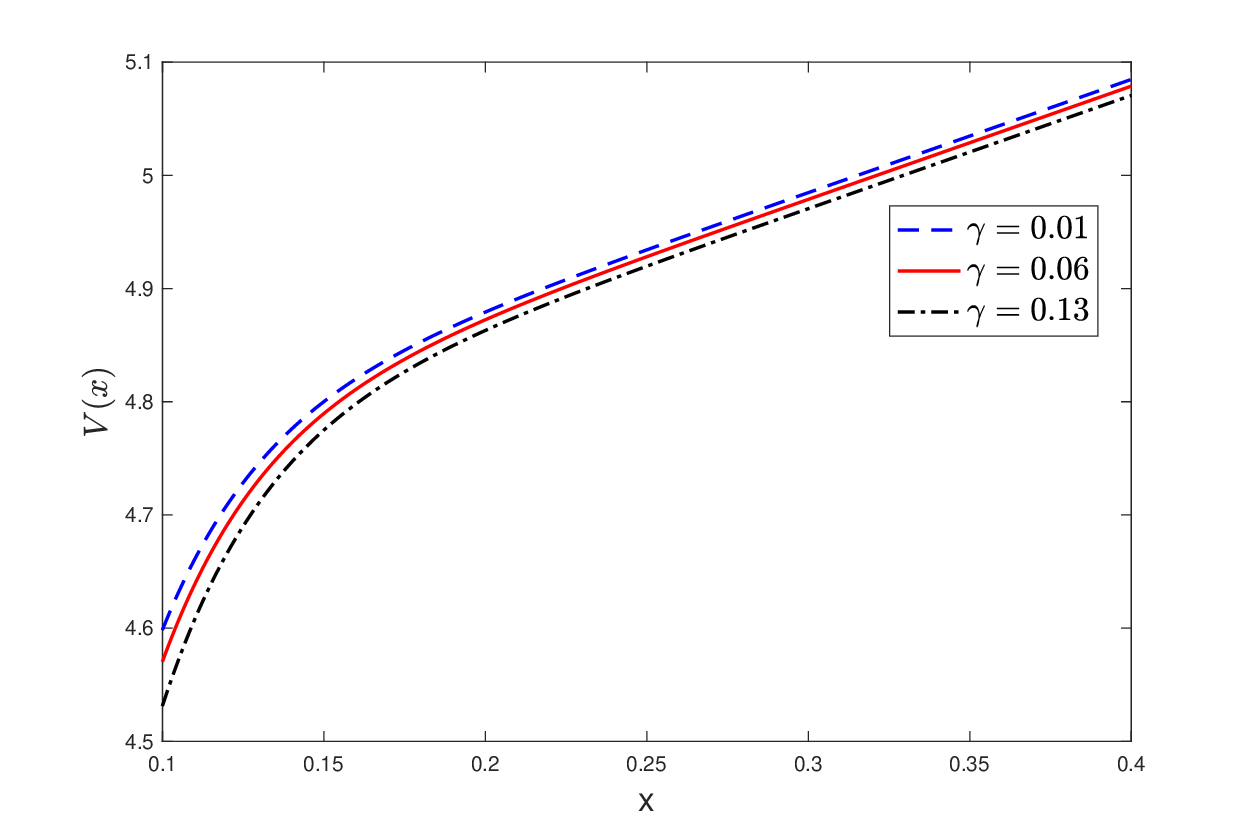}
	\end{subfigure}
\vspace{-4ex}
	\caption{The value function $V(x)$ (left) and its zoom-in for $x \in (0.1, 0.4)$ (right)}
	\label{fig2}
\end{figure}

Next, we compute the equilibrium value function $V$ in \eqref{eq:Vx} for three different risk aversion levels, $\gamma = 0.01, 0.06$, and $0.13$, and plot their graphs as a function of $x$ in Figure \ref{fig2}. The left panel verifies the strict concavity of $V$ over $[0, \xt)$, while the right panel shows that $V$, viewed as a function of $\gamma$, is decreasing.   We expect $V$ to decrease with respect to $\gam$ because of the form of the objective function $J$ in \eqref{eq:J}, and it is satisfying to see our intuition born out in Figure \ref{fig2}.

Recall that the equilibrium consumption strategy in Kronborg and Steffensen \cite{KS2015} is of a bang-bang type, depending solely on the order between the risk-free rate and the discount rate $\rho$, and it is independent of the state process. This motivates us to study the impact of $\rho$ on our results. In particular, we study how $\rho$ affects the upper bound $\bar{\gam}$ (Theorem \ref{thm:onebarrier} requires $\gam \le \bar{\gam}$) and the barrier $\xt := \xt (\rho)$ for a given small risk aversion $\gam = 0.1396$. We first observe a technical result that the upper bound on risk aversion, $\bar{\gam}$, in Theorem \ref{thm:onebarrier} increases with respect to $\rho$. The right panel shows that the (unique) barrier $\tilde{x}$ decreases as $\rho$ increases,  indicating that when the discount rate is higher (that is, the manager is less patient), larger dividends are paid out earlier.  

\begin{figure}[h]
	\begin{subfigure}{0.5\textwidth}
		\includegraphics[width = 8cm, height = 5cm, trim = 1cm 0cm 1cm 0.8cm, clip = true]{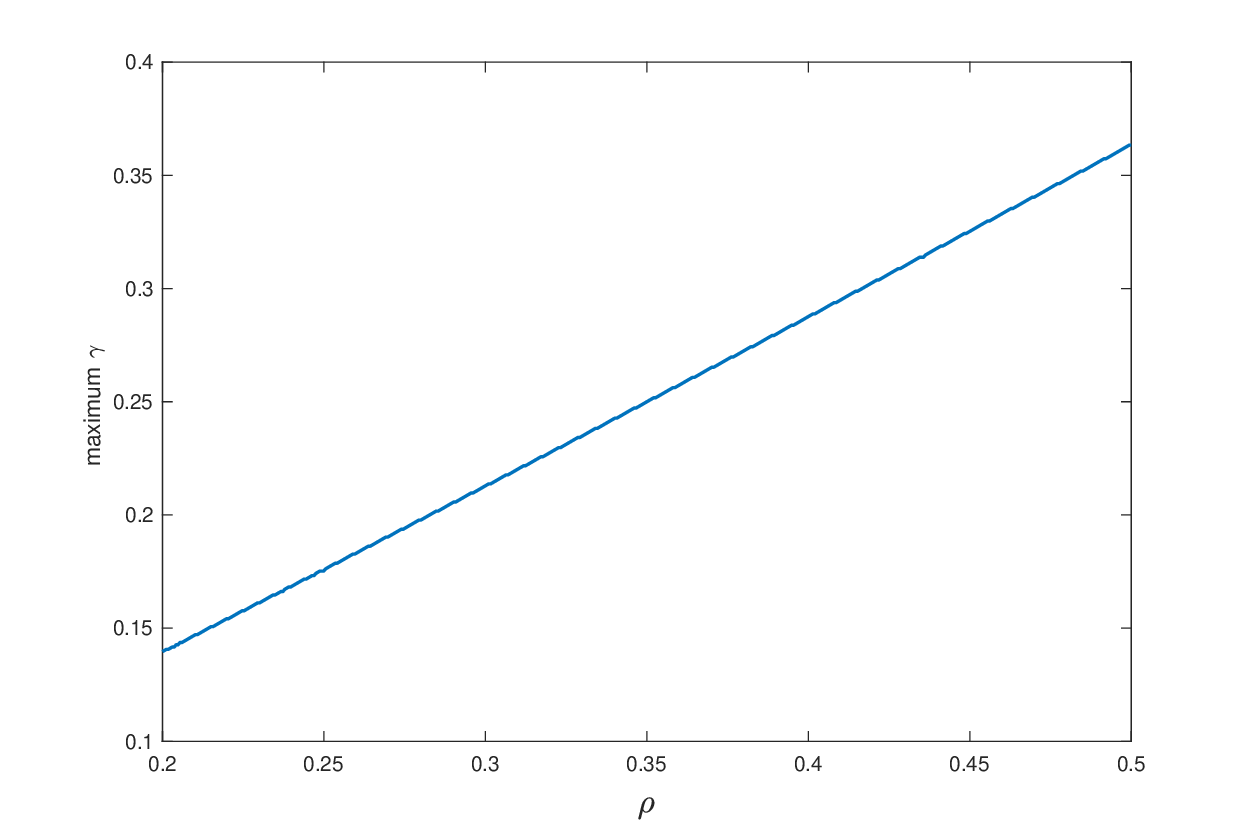}
	\end{subfigure}
	\begin{subfigure}{0.5\textwidth}
		\includegraphics[width = 8cm, height = 5cm, trim = 1cm 0cm 1cm 0.8cm, clip = true]{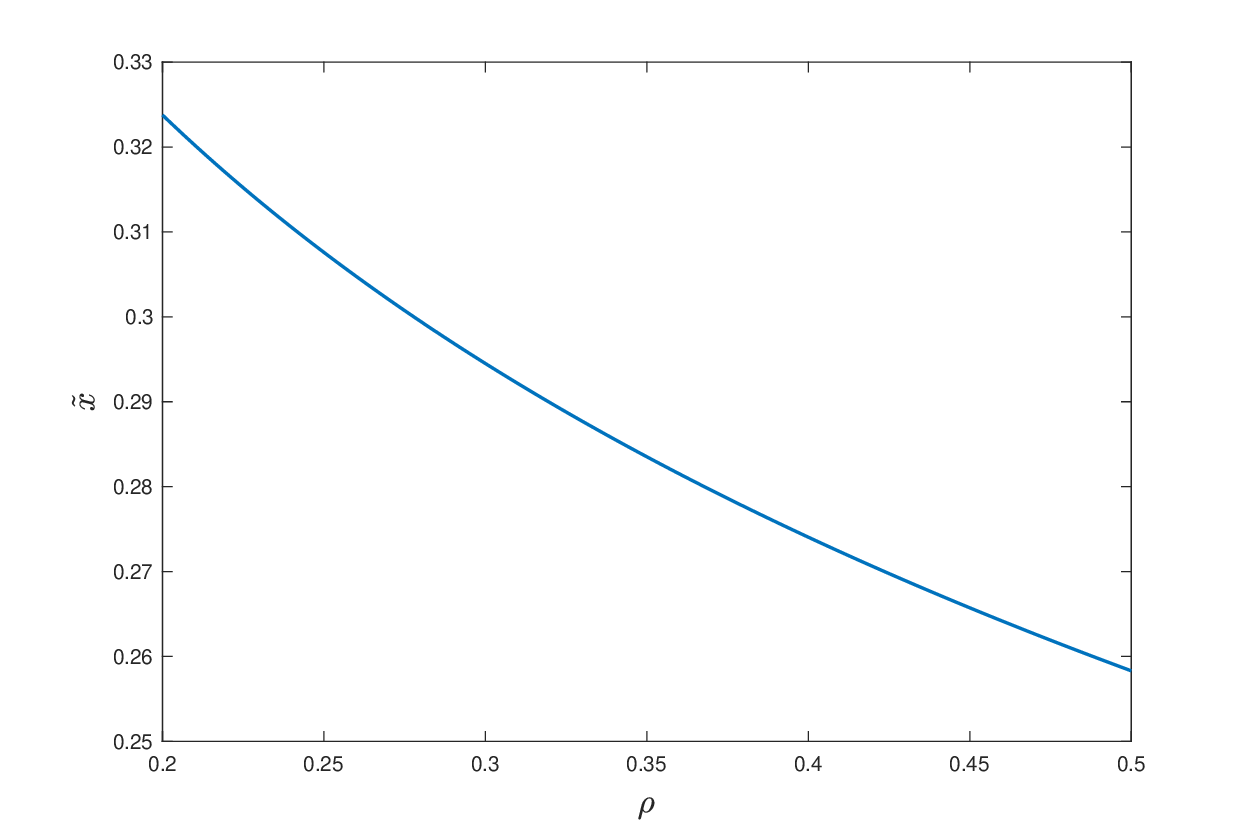}
	\end{subfigure}
	\vspace{-4ex}
	\caption{Impact of the discount rate $\rho$ on the maximum allowed risk aversion $\bar{\gam}$ (left) and the barrier $\xt$ under $\gam = 0.1396$ (right)}
	\label{fig:rho}
\end{figure}

\begin{figure}[h]
	\begin{subfigure}{0.5\textwidth}
		\includegraphics[width = 9cm, height = 6cm]{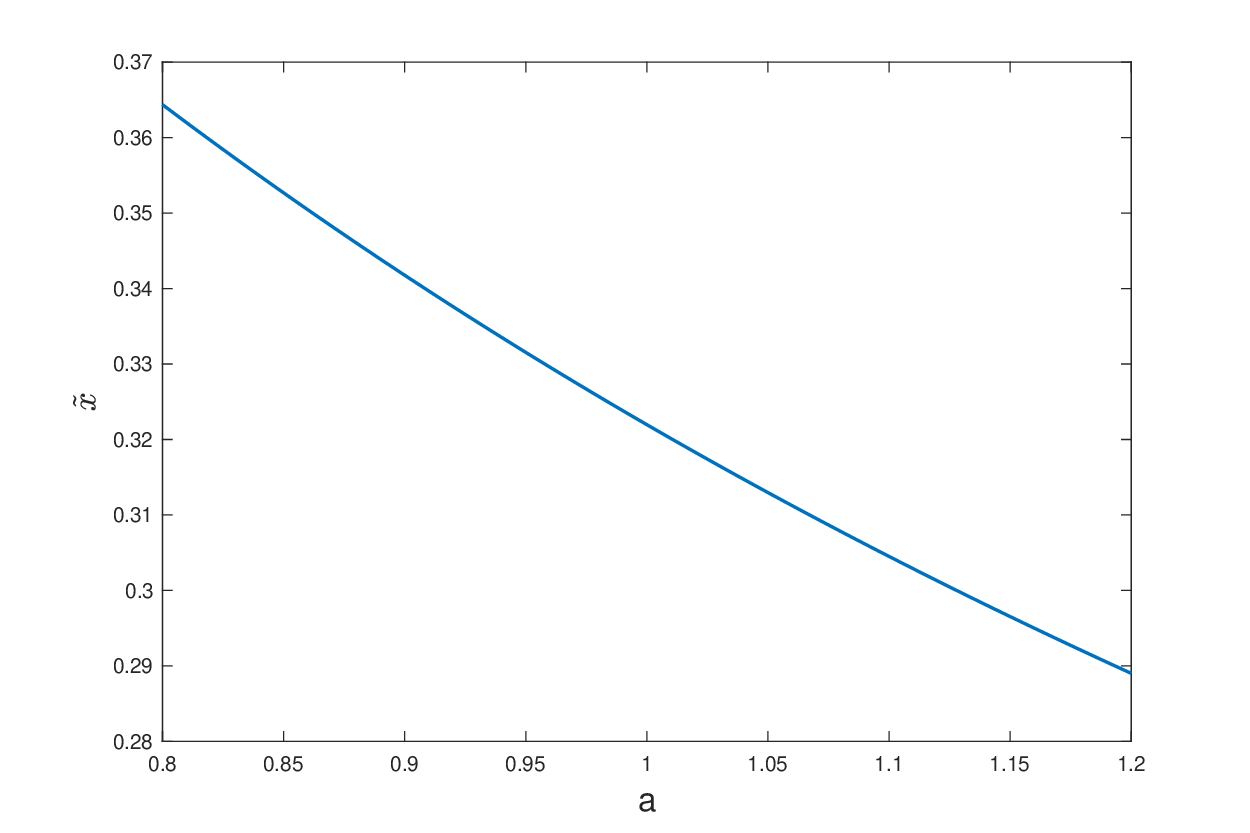}
	\end{subfigure}
	\begin{subfigure}{0.5\textwidth}
		\includegraphics[width = 9cm, height = 6cm]{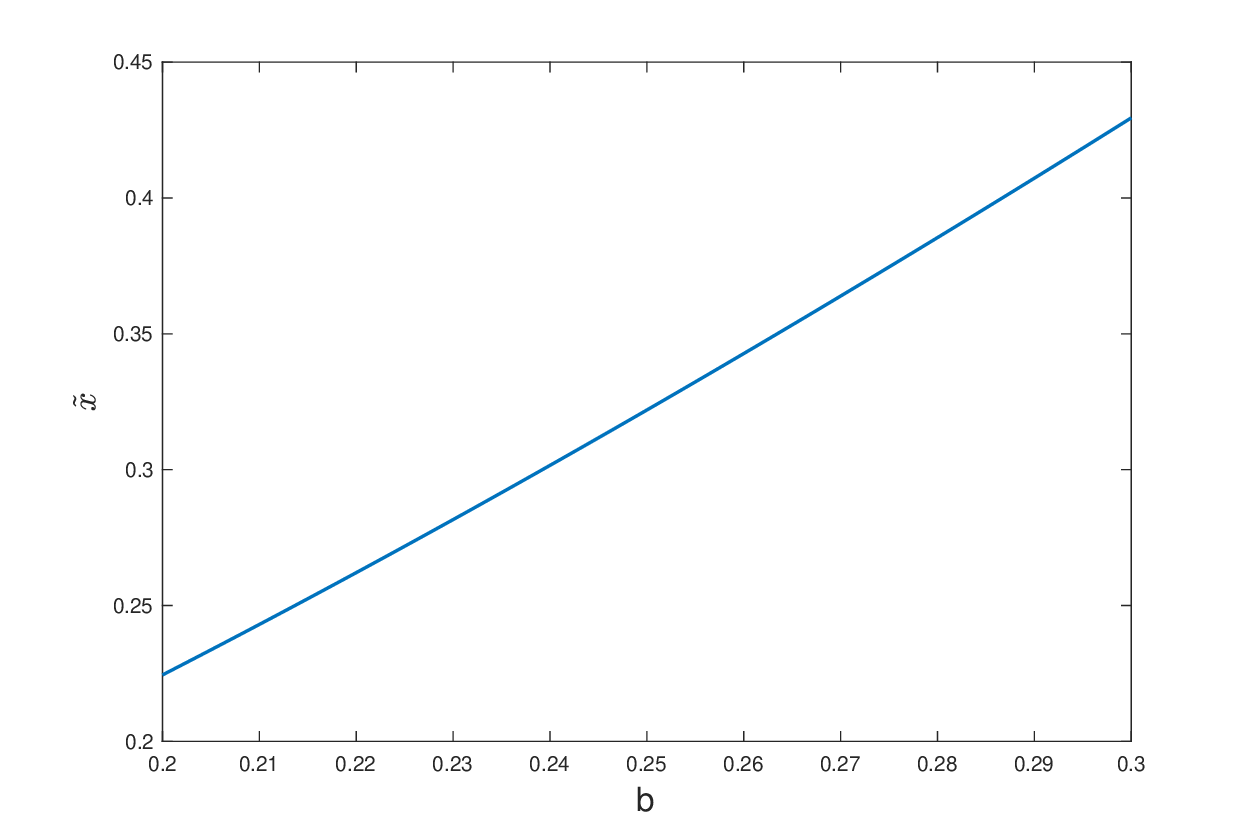}
	\end{subfigure}
	\caption{Impact of the drift $a$ (left) and volatility $b$ (right) on the barrier $\xt$ when $\gam = 0.11$}
	\label{fig:ab}
\end{figure}

We proceed to study how the model parameters $a$ and $b$ in the diffusion surplus influence the barrier $\tilde{x}$ given by \eqref{eq:xt}. Again, we focus on the cases of small risk aversion $\gam$ in which Theorem \ref{thm:onebarrier} applies, and note that the admissible range of $\gam$ is implicitly determined by the model parameters. Here, we present results for $\gam = 0.11$ in Figure \ref{fig:ab}, under which all conditions of Theorem \ref{thm:onebarrier} are satisfied over the range of $a$ and $b$ considered. Together with  Figure \ref{fig1} (left panel) and Figure \ref{fig:rho} (left panel), these results indicate that the barrier $\tilde{x}$ varies continuously with respect to $\gam$, $\rho$, $a$ and $b$, as long as the assumptions of Theorem \ref{thm:onebarrier}  hold.

\begin{figure}[h]
	\begin{subfigure}{0.5\textwidth}
		\includegraphics[width = 9cm, height = 6cm]{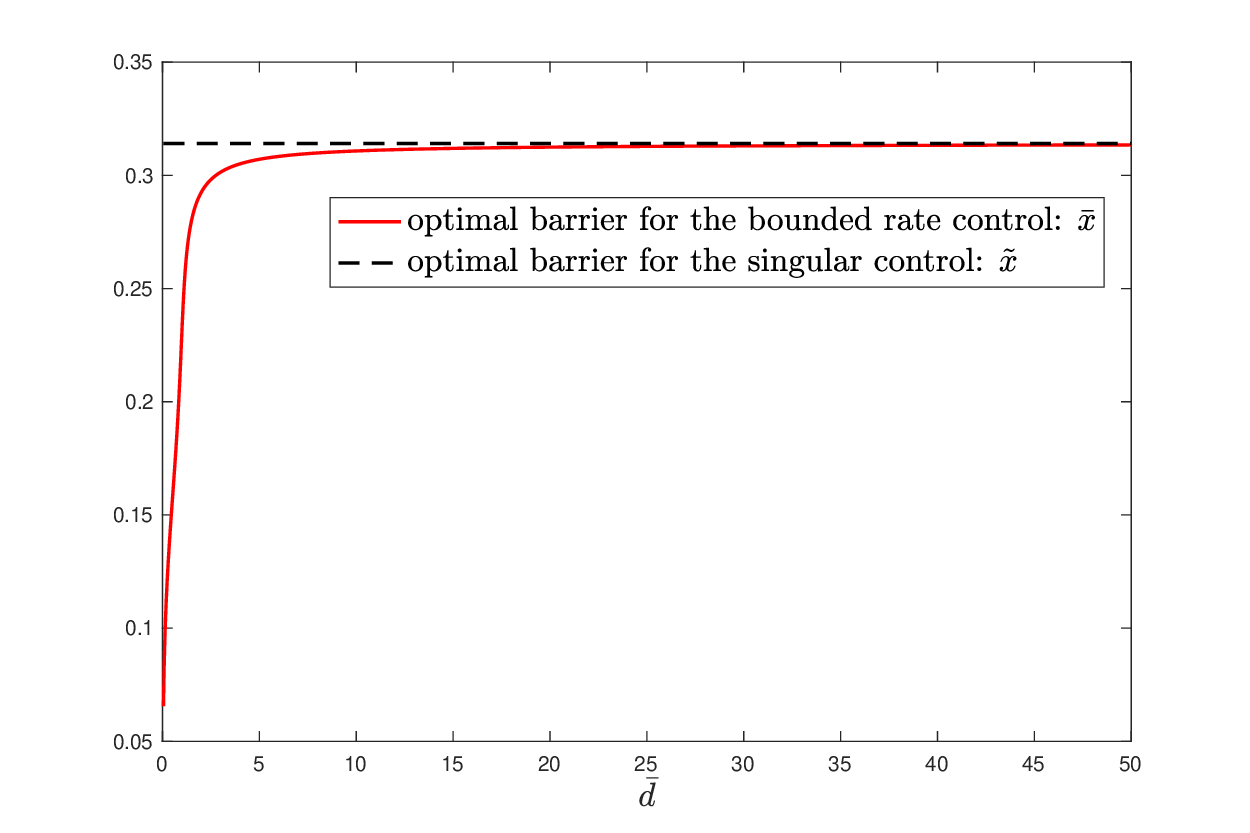}
	\end{subfigure}
	\begin{subfigure}{0.5\textwidth}
		\includegraphics[width = 9cm, height = 6cm]{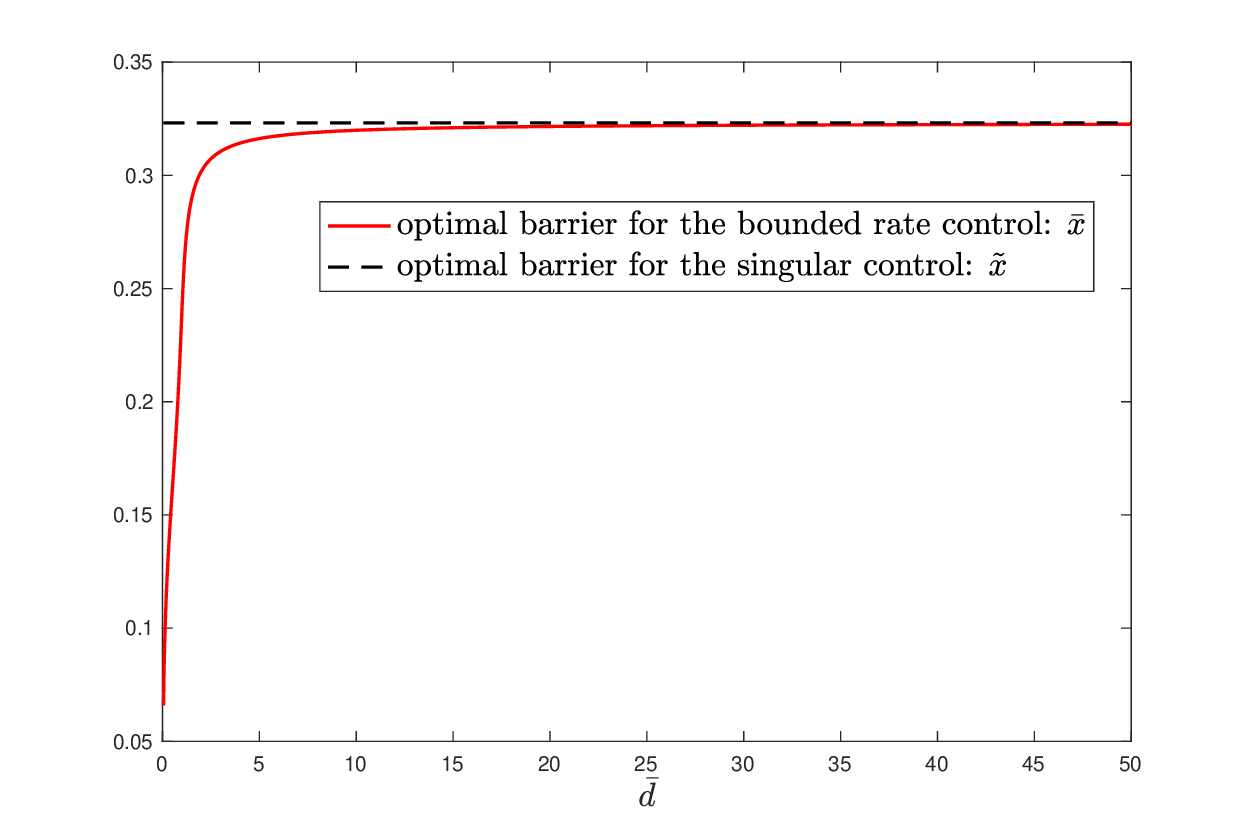}
	\end{subfigure}
	\caption{The barriers for the classical control and singular control problems when $\gam = 0$ (left) and $\gam = 0.13$ (right)}
	\label{fig:rate_db}
\end{figure}

Cao et al.\ \cite{CLYZ2025} study a similar MV dividend problem under the classical control framework and seek an equilibrium dividend \emph{rate} strategy, subject to a maximum payout rate $\bar{d}>0$. For a fixed (feedback) dividend rate strategy $\pd$, the cumulative dividend $D_t$ is given by  
 	\begin{align*}
 		D_t = \int_0^{t\wedge \tau} \pd(X_s, s) \drm s.
 	\end{align*}
They show that, for sufficiently small risk aversion $\gamma$ (along with conditions ensuring the uniqueness of a positive solution $\bar{x}$ to a nonlinear equation), a barrier strategy $\pd^*(x) = \bar{d}\,\id_{{x > \bar{x}}}$ is an equilibrium strategy; this result is parallel to ours in Theorem \ref{thm:onebarrier}, in which $\xt$ is the unique barrier. To examine the connections between two different frameworks (singular control versus classical control), we set the same parameters for $a$, $b$, and $\rho$ as above and compute the two barriers $\xt$ and $\bar{x}$ (for the latter, we consider $\bar{d} \in [0, 50]$). We plot their graphs when $\gamma = 0$ (left panel) and $\gamma = 0.13$ (right panel) in Figure \ref{fig:rate_db}. We observe that, as the maximum dividend rate $\bar{d}$ increases, the corresponding barrier $\bar{x}$ in Cao et al.\ \cite{CLYZ2025} converges to $\xt$ of this paper. In fact, when $\gamma = 0$, Jeanblanc-Picqu\'e and Shiryaev \cite{JS1995} prove that as $\bar{d} \to \infty$, the \emph{optimal} barrier of the bounded-rate problem converges to that of the singular control problem. Our numerical results suggest that this convergence holds for \emph{equilibrium strategies} when $\gamma$ is small.

\section{Conclusions}
\label{sec:con}

In this paper, we studied a novel singular control, time-inconsistent dividend problem, and the objective is to optimize the MV criterion of the integral of all discounted dividends paid until ruin time, an endogenous stopping time. 
We proved a new verification theorem that characterizes equilibrium dividend strategies and their corresponding value functions to this problem.
We, then, used the verification theorem to prove two results in which we obtain equilibrium dividend strategies (semi-)explicitly: one for large values of risk aversion $\gam$ (specifically, $\gam \ge \frac{2a}{b^2}$), and one for small values of $\gam$ (namely, $\gam < \eps \le \frac{2a}{b^2}$, subject to Condition (2) of Theorem \ref{thm:onebarrier}, in which $\eps$ depends upon the parameters of the model). Numerical experiments show that the maximum $\gam$ satisfying both conditions of Theorem \ref{thm:onebarrier}, denoted by $\bar{\gam}$,  is strictly less than $\frac{2a}{b^2}$. Thus, finding equilibrium dividend strategies when $\bar{\gam} < \gam < \frac{2a}{b^2}$ remains an open question. 

For future work, one direction is to allow investment or capital injection strategies, in addition to dividend control, in the model (see Albrecher and Thonhauser \cite{AT2009}). Recall that we study dividend control problems up to the ruin time in this work; however, there are alternative definitions of ``ruin'' (see Section 5 in Avanzi \cite{A2009}), and it will be interesting to revisit our problem under such alternative definitions. 
In this paper, we chose the notion of weak equilibrium (see Bj\"ork and Murgoci \cite{BM2010}), and several recent papers pointed out its potential drawback and proposed different notions of equilibrium, such as strong equilibrium (see, for instance, Bayraktar et al.\ \cite{BZZ2021}, Bayraktar et al.\ \cite{BWZ2023}, and He and Jiang \cite{HJ2021}). To the best of our knowledge,  time-inconsistent singular control problems under the notion of strong equilibrium have not been studied before. Our numerical experiments show that the barrier $\xt$ varies continuously with respect to the model parameters (see Figures \ref{fig1}, \ref{fig:rho}, and \ref{fig:ab}), when the conditions of Theorem \ref{thm:onebarrier} hold; however, an \emph{analytical} study of the stability of equilibria (as in Bayraktar et al.\ \cite{BWZ2023b}) remains an open question, and we leave it for future research.

\vspace{2ex}
\noi
\textbf{Acknowledgments.}  We thank the corresponding editor, Professor Erhan Bayraktar, and anonymous associate editor and reviewers for their valuable comments. The first and second authors acknowledge the financial support from the Natural Sciences and Engineering Research Council of Canada, grants 05061 and 04958, respectively. The third author thanks the Cecil J. and Ethel M. Nesbitt Professorship for the financial support of her research.

\appendix

\small 
\setlength{\bibsep}{2pt}

\end{document}